%% file: HNA_BEM_transmission_Groth_Hewett_Langdon_WM_revision.tex
\documentclass[12pt]{article}

\usepackage{amsmath}
\usepackage{amsfonts}
\usepackage{amsthm}
\usepackage{amssymb}
\usepackage{enumerate,verbatim}
\usepackage{cases}
\usepackage{multirow}
\usepackage{subfigure}
\usepackage{pgf,tikz}
\usepackage{pgfplots}
\pgfplotsset{compat=1.4}
\pgfplotsset{every axis plot/.append style={line width=0.8pt}}
\usepackage{graphicx}
\usepackage{appendix}
\usetikzlibrary{calc}
\usepackage{color}
\usepackage[hide]{ed} 
\usetikzlibrary{arrows,positioning}
\definecolor{dhcol}{rgb}{0,0.5,0}
\providecommand{\dhc}[1]{{\color{dhcol}#1}}
\definecolor{sgcol}{rgb}{0,0,0.5}
\providecommand{\sgc}[1]{{\color{sgcol}#1}}
\definecolor{slcol}{rgb}{0.5,0,0}

\newcommand{\dhednote}[1]{\ednote{\dhc{{DH:} #1}}} 
\newcommand{\sgednote}[1]{\ednote{\sgc{{SG:} #1}}}

\usetikzlibrary{%
    decorations.pathreplacing,%
    decorations.pathmorphing%
}

\begin{document}

\input{macros}
\newcommand{\HH}{\mathbf{H}}
\newcommand{\EE}{\mathbf{E}}
\newcommand{\kk}{\mathbf{k}}
\newcommand{\be}{\mathbf{e}}
\newcommand{\bof}{\mathbf{f}}
\newcommand{\bv}{\mathbf{v}}
\newcommand{\bw}{\mathbf{w}}
\newcommand{\bX}{\mathbf{X}}
\newcommand{\bt}{\mathbf{t}}
\newcommand{\sd}{\mbox{d}}
\newcommand{\dr}{\bd}
\newcommand{\di}{\be}
\newcommand{\dri}{\dr^i}
\newcommand{\dii}{\di^i}
\newcommand{\drr}{\dr^r}
\newcommand{\dir}{\di^r}
\newcommand{\drt}{\dr^t}
\newcommand{\dit}{\di^t}
\newcommand{\thetari}{\theta^i}
\newcommand{\thetaii}{\phi^i}
\newcommand{\thetarr}{\theta^r}
\newcommand{\thetair}{\phi^r}
\newcommand{\thetart}{\theta^t}
\newcommand{\thetait}{\phi^t}
\newcommand{\NN}{D}
\newcommand{\KK}{E}
\newcommand{\Ns}{{\tilde{\NN}_i}}
\newcommand{\Ks}{{\tilde{\KK}_i}}
\newcommand{\DOFone}{\mathrm{\#DOFper}\lambda_1}
\newcommand{\DOFtwo}{\mathrm{\#DOFper}\lambda_2}
\newcommand{\sumpm}[1]{\{\!\!\{#1\}\!\!\}}

\title{A hybrid numerical-asymptotic boundary element method for high frequency scattering by penetrable convex polygons}
\author{S. P. Groth$^{\text{a}}$\footnotemark[1]\footnote{Present address: Department of Electrical Engineering and Computer Science, Massachusetts Institute of Technology, 77 Massachusetts Ave., Cambridge, Massachusetts 02139, United States of America.}\ \footnotemark[2]\footnote{Corresponding author. Email address: samgroth@mit.edu},
 D. P. Hewett$^{\text{b}}$ and S. Langdon$^{\text{a}}$\\[2pt]
$^{\text{a}}$Department of Mathematics and Statistics, University of Reading, \\[4pt]
 Whiteknights PO Box 220, Reading, RG6 6AX, United Kingdom.\\[6pt]
$^{\text{b}}$Department of Mathematics, University College London,\\[2pt]
Gower Street, London, WC1E 6BT, United Kingdom.
 }
\maketitle

\begin{abstract}
We present a novel hybrid numerical-asymptotic boundary element method for high frequency acoustic and electromagnetic scattering by 
penetrable (dielectric) convex polygons. 
Our method is based on a standard reformulation of the associated transmission boundary value problem as a direct boundary integral equation for the unknown Cauchy data, but with a nonstandard numerical discretization which efficiently captures the high frequency oscillatory behaviour. 
The Cauchy data is represented as a sum of the classical geometrical optics approximation, computed by a beam tracing algorithm, plus a contribution due to diffraction, computed by a Galerkin boundary element method using oscillatory basis functions chosen according to the principles of the Geometrical Theory of Diffraction. 
We demonstrate with a range of numerical experiments that our boundary element method can achieve a fixed accuracy of approximation using only a relatively small, frequency-independent number of degrees of freedom. 
Moreover, for the scattering scenarios we consider, 
the inclusion of the diffraction term provides an order of magnitude improvement in accuracy over the geometrical optics approximation alone. 
\end{abstract}


\section{Introduction}
\label{sec:intro}
Scattering of time-harmonic acoustic or electromagnetic waves by penetrable (i.e., dielectric) scatterers arises in numerous physical applications, for example  
the scattering of light by atmospheric ice crystals, important in determining the earth's radiation balance in global climate models~\cite{Baran:2009}. 
When the penetrable scatterer and the exterior domain are both homogeneous,
a standard computational 
approach to solving the 
associated transmission boundary value problem, at least for low/moderate frequencies, is the boundary element method (BEM) (often called the Method Of Moments in the electromagnetic community), which involves the numerical solution of a system of boundary integral equations (BIEs) on the scatterer boundary \cite{sauterschwab,smigaj2012solving,groth2015boundary}. 

However, conventional BEM implementations based on piecewise polynomial approximation spaces are computationally infeasible in the \emph{high frequency} regime where the scatterer is large relative to the smallest wavelength in the problem, since they suffer from the well-known limitation (suffered by all conventional numerical approaches including finite element, finite difference and spectral methods) that a fixed number of degrees of freedom is required per wavelength in order to accurately represent the oscillatory solution \cite{marburg2008discretization}. At high frequencies the wave field is more naturally described by asymptotic theories such as Geometrical Optics (GO) and the Geometrical Theory of Diffraction (GTD) \cite{Keller,Kinber,james1986geometrical}, which represent it in a completely different way, as a superposition of ray fields. In principle, these asymptotic approximations have frequency-independent computational cost. However, they are accurate only at sufficiently high frequencies, and require the solution of certain canonical problems, which are not always solvable in closed form.

In this paper we present a \emph{hybrid numerical-asymptotic} (HNA) BEM for two-dimensional scattering by penetrable convex polygons. 
In general terms, the HNA 
approach aims to develop numerical methods 
which are computationally feasible over the whole frequency range, by explicitly building features of the high frequency solution behaviour into the numerical approximation space. 
Specifically, rather than using conventional piecewise polynomial BEM basis functions, one instead uses oscillatory basis functions, whose oscillations match those of the GO and GTD ray fields. 
The HNA approach 
(reviewed in \cite{acta}) 
has proven very effective for a number of problems involving impenetrable scatterers, including convex \cite{CWL:07,hewett2013high,CWLM}, nonconvex \cite{chandler2012high,hewett:shadow} and curvilinear \cite{LaMoCh:10} polygons, two-dimensional colinear screens \cite{hewett2014frequency}, smooth convex two-dimensional \cite{dominguez2007hybrid,bruno2004prescribed,EcOz:17,EcUr:16} and three-dimensional \cite{GH11} scatterers, and three-dimensional rectangular screens \cite{Hargreaves}. 
All of the HNA methods mentioned above 
achieve fixed accuracy with a frequency-independent (or only very modestly growing as the frequency increases) number of degrees of freedom. 
Moreover, many 
(in particular \cite{CWL:07,hewett2013high,chandler2012high,hewett:shadow,hewett2014frequency,EcOz:17,EcUr:16})
are fully convergent and supported by rigorous numerical analysis. 

We believe that the method presented in this paper is the first HNA method for any problem involving a penetrable scatterer.  
Application of the HNA methodology to penetrable scatterers is nontrivial, primarily because the high frequency asymptotic behaviour is considerably more complicated and less well-understood than in the corresponding impenetrable case. In addition to the phenomena of reflection and diffraction that pertain in the impenetrable case, one also has to take into account the \emph{refraction} of waves into the interior of the scatterer. This generates a complicated superposition of internally reflected and diffracted wave fields. Furthermore, the exact nature of these wave fields is not fully understood, because many of the relevant GTD canonical scattering problems remain unsolved in closed form. 
In particular, the key open problem of relevance for scattering by penetrable polygons is diffraction by a penetrable wedge, which has not yet been solved in full generality, despite numerous attempts (see, e.g., \cite{MePeSpTe:94,Rawlins,budaev1999rigorous,antipov2007diffraction}).
Nonetheless, in spite of these difficulties it was demonstrated in our earlier publication \cite{groth2015hybrid}
that effective HNA approximation spaces can indeed be developed for penetrable scatterers, by adopting appropriate generalizations of the methods already developed for the impenetrable case. 

The two key general ideas explored in \cite{groth2015hybrid} are that:
\begin{enumerate}[(i)]
\item In order to design HNA approximation spaces one needs only \emph{partial }information about the GTD solution (specifically phases, and the location of shadow boundaries), not the full GTD solution. In particular, the lack of analytical expressions for corner diffraction coefficients for penetrable scatterers is not a barrier to the development of HNA methods;
\item Even if the high frequency oscillatory behaviour is too complicated to be captured completely by an HNA method (because there are too many phases to consider), incorporating just a small number of oscillatory basis functions capturing the leading order asymptotic behaviour might lead to a significant improvement in performance over conventional purely numerical or purely asymptotic methods.
\end{enumerate}
Specifically, for the case of penetrable convex polygons, the approximation space proposed in \cite{groth2015hybrid} represents the unknown Cauchy data in the BIE formulation as the GO field (which can be computed analytically using a beam-tracing method - cf.\ \S\ref{sec:comp_v_go} below) plus an HNA component designed to capture the primary, secondary, and all higher-order corner-diffracted fields. This HNA component takes the form of a sum of terms $\bv_1^j\re^{\ri k_1r_j}$ and $\bv_2^j\re^{\ri k_2r_j}$ (see equation \rf{eqn:ansatztransmission} below), where $k_1$ and $k_2$ are the exterior and interior wavenumbers, $r_j$ is Euclidean distance from the $j$th corner of the polygon, and $\bv_1^j$ and $\bv_2^j$ are slowly-varying amplitudes, which are represented by piecewise polynomials on appropriately graded meshes. 
A BEM implementation was not presented in \cite{groth2015hybrid}; instead  the performance of the proposed approximation space was investigated by projecting a highly accurate fully numerical reference solution onto the approximation space and computing the resulting approximation error.

In the current paper we present a Galerkin BEM implementation of a refined version of the HNA approximation space developed in \cite{groth2015hybrid}, demonstrating that it can be used as the basis of an effective numerical method. Regarding the approximation space itself, compared to \cite{groth2015hybrid} we have improved: (i) the robustness of the beam-tracing algorithm used to compute the classical GO approximation (see \S\ref{sss:sign}), leading to improved accuracy for certain incident wave directions; and (ii) the efficiency of our meshing strategy for treating beam boundary discontinuities (see \S\ref{sec:BeamBoundaries} and \S\ref{sec:Vd}), significantly reducing the number of degrees of freedom required to approximate the diffracted wave fields. 

The state of the art for high frequency electromagnetic transmission problems in atmospheric physics appears to be the physical-geometric optics hybrid (PGOH) method of Bi et al.\ \cite{bi2013physical,Bi:scat}. The PGOH method is an extension to penetrable scatterers of the classical Kirchhoff approximation, in which one takes a boundary integral representation of the solution in the propagation domain (Green's representation theorem in acoustics; the Stratton-Chu formula in electromagnetics) and replaces the unknown boundary data by its GO approximation (including contributions not just from the incident wave, but also from internally reflected and refracted rays).
Our method improves on the PGOH approximation by adding in numerically computed diffracted fields within the integral representation. 
We shall demonstrate, by a range of numerical experiments that, for all of the examples we consider, our method achieves an order of magnitude improvement in accuracy over the PGOH approach, with only a modest (and frequency-independent) number of degrees of freedom required to compute the additional diffracted contribution. 

The outline of the paper is as follows. We begin in~\S\ref{sec:problem} by stating the scattering problem and its reformulation as a BIE.  In \S\ref{sec:HNA} we describe the construction of our HNA approximation space. 
In \S\ref{sec:Numerical} we describe the implementation of our approximation space as a Galerkin BEM. We present numerical results demonstrating its efficiency and accuracy for a range of examples with different (convex polygonal) geometries and different physical parameters, including variations in contrast, absorption, and incident angle. 
In \S\ref{sec:Discussion} we offer some conclusions and discuss potential future improvements and generalizations of the method.

\section{Problem statement}
\label{sec:problem}
We consider the two-dimensional scattering of a time-harmonic incident plane wave by a penetrable convex polygon. 
Let $\Omega_2\subset\R^2$ denote the interior of the polygon (a convex bounded open set), and let $\Omega_1:=\mathbb{R}^2\backslash\overline{\Omega_2}$ denote the exterior unbounded domain. Let $\Gamma=\Gamma_1\cup\Gamma_2\cup\ldots\Gamma_{n_s}$ denote the boundary of the polygon, where $n_s$ is the number of sides and $\Gamma_j$, $j=1,\ldots,n_s$, are the sides of the polygon, labelled in an anti-clockwise direction. 
The corners of the polygon are similarly labelled $\bP_1, \ldots,\bP_{n_s}$, with $\Gamma_j$, $j=1,\ldots,n_s$, being the side between the corners $\bP_j$ and $\bP_{j+1}$ (with the convention that $\bP_{n_s+1}\equiv \bP_1$); for an example see Figure \ref{fig:triangle_setup}(a). 
Let $\bn$ denote the outward unit normal to $\Gamma$. Let $k_1$ denote the wavenumber in 
$\Omega_1$ and $k_2$ 
the wavenumber in 
$\Omega_2$. We shall assume throughout that
\begin{align}
\label{k1k2Assumptions}
k_1>0,\, \real{k_2}>0, \text{ and } \im{k_2}\geq 0;
\end{align}
when $\im{k_2}>0$ the scatterer is partially absorbing. For convenience we introduce the notation
\begin{align*}
\label{}
\mu = k_2/k_1
\end{align*}
for the (possibly complex) refractive index, and $\lambda_1=2\pi/k_1$, $\lambda_2=2\pi/\real{k_2}=\lambda_1/\real{\mu}$ for the exterior and interior wavelengths.

The boundary value problem (BVP) we wish to solve is: given $\alpha\in\C\setminus \{0\}$ satisfying
\begin{align}
\label{alphaAssumptions}
\im{\alpha}\leq 0,\,\im{\alpha k_2^2}\geq0,
\end{align}
and an incident field
\begin{equation}
   u^i(\bx):=\re^{\ri k_1\bd^i \cdot \bx}, \qquad \bx
   \in\R^2,
   \label{eqn:ui}
\end{equation}
with $\bd^i\in\R^2$ a unit direction vector,  
determine the total field $u_1$ in $\Omega_1$ and $u_2$ in $\Omega_2$ such that 
\begin{align}
\Delta u_1 + k_1^2u_1 = 0, & \qquad \mbox{in} \ \Omega_1, \label{helm_ext}\\ 
\Delta u_2 + k_2^2u_2 = 0, &\qquad \mbox{in} \ \Omega_2,\label{helm_int}\\
u_1=u_2\ \mbox{and}\ \frac{\partial u_1}{\partial\bn}=\alpha\frac{\partial u_2}{\partial\bn}, &\qquad  \mbox{on}\ \Gamma, \label{bcs}
\end{align}
with the scattered field $u^s:=u_1-u^i$ satisfying the Sommerfeld radiation condition
\begin{equation}
    \frac{\partial u^s}{\partial r}(\bx) - \ri k_1 u^s(\bx) = o(r^{-1/2}), \quad \mbox{as }r:=|\bx|\rightarrow\infty.
\label{SRC}
\end{equation}

The scalar BVP \rf{eqn:ui}--\rf{SRC} is a model for both acoustic and electromagnetic scattering problems. In the electromagnetic case, scattering by a (possibly partially absorbing) dielectric polygon is modelled by two BVPs of the form \rf{eqn:ui}--\rf{SRC}, one for the out-of-plane electric field and one for the out-of-plane magnetic field. The standard transmission boundary conditions for Maxwell's equations (see, e.g., \cite[\S1.1]{Born}) imply that for the electric field the appropriate choice of $\alpha$ in \rf{bcs} is $\alpha=1$; for the magnetic field it is $\alpha=1/\mu^2=(k_1/k_2)^2$. Assuming \rf{k1k2Assumptions}, both of these choices of $\alpha$ satisfy \rf{alphaAssumptions}. 

Under the assumptions \rf{k1k2Assumptions} and \rf{alphaAssumptions} it is well known that the BVP \rf{eqn:ui}--\rf{SRC} is uniquely solvable with $u_1\in C^2(\Omega_1)\cap H^{1}_{\rm loc}(\Omega_1)$ and $u_2\in C^2(\Omega_2)\cap H^{1}(\Omega_2)$ 
(see, e.g., \cite[Proposition~2.1 and Corollary~3.4]{LaRaSa:09}, which follows from results in \cite{costabel,torres1993helmholtz}, and also the related result of \cite[Corollary~8.5]{marmolejo2012transmission}). We mention also \cite{MoSp:17}, where a particularly simple proof of well-posedness is given for the case $k_1,k_2>0$, along with wavenumber-explicit stability bounds.

Our BEM for solving \rf{eqn:ui}--\rf{SRC} is based on a direct BIE formulation. We first observe that if $u_1$ and $u_2$ satisfy \rf{eqn:ui}--\rf{SRC}, then a form of Green's representation theorem holds; i.e., the fields $u_1$ and $u_2$ in $\Omega_1$ and $\Omega_2$ can be represented in terms of their Dirichlet and Neumann traces on $\Gamma$.
By applying the transmission conditions \rf{bcs} we can work with the boundary traces of just one of the fields; we choose to work with $u_1$, denoting its Dirichlet and Neumann traces on $\Gamma$ by $u_1$ and $\pdonetext{u_1}{\bn}$ respectively. 
Then 
\cite[Theorems~2.20 and~2.21]{acta}
\begin{alignat}{3}
u_1(\bx) &= u^i(\bx) - \cS_1 \pdone{u_1}{\bn}(\bx)+\cD_1 u_1(\bx), &\qquad \bx\in\Omega_1, \label{eqn:u1} \\ 
u_2(\bx) &= \frac{1}{\alpha}\cS_2 \pdone{u_1}{\bn}(\bx)-\cD_2 u_1(\bx) &\qquad \bx\in\Omega_2, \label{eqn:u2}
\end{alignat}
where $\cS_j:H^{-1/2}(\Gamma)\to H^1_{\rm loc}(\R^2\setminus \Gamma)$ and $\cD_j:H^{1/2}(\Gamma)\to H^1_{\rm loc}(\R^2\setminus \Gamma)$, $j=1,2$, are the standard single- and double-layer potentials, defined for sufficiently smooth arguments (for a detailed discussion of the properties of $\cS_j$ and $\cD_j$ on Lipschitz domains see \cite{acta}) by the integral formulas
\begin{align*}
\label{}
\cS_j\phi(\bx)=\int_\Gamma \Phi_j(\bx,\by)\phi(\by) \, \rd s(\by),
\qquad
\cD_j\phi(\bx)=\int_\Gamma \pdone{\Phi_j(\bx,\by)}{\bn(y)}\phi(\by) \, \rd s(\by).
\end{align*}
Here $\Phi_j(\bx,\by):=(\ri/4)H_0^{(1)}(k_j|\bx-\by|)$, $j=1,2$, are the fundamental solutions of the Helmholtz equations (\ref{helm_ext}) and (\ref{helm_int}), respectively, with $H_{0}^{(1)}$ denoting the Hankel function of the first kind of order zero. 
For later reference we note that the far-field behaviour of the scattered field is\begin{align*}
u^s(\bx) \sim \frac{\re^{\ri\pi/4}}{2\sqrt{2\pi}}\frac{\re^{\ri k_1r}}{\sqrt{k_1r}}F(\hat{\bx}),\quad \mbox{as}\ r:=|\bx|\rightarrow\infty,
\end{align*}
where $\hat{\bx}:=\bx/|\bx|\in\mathbb{S}^1$, the unit circle, and the \textit{far-field pattern} $F:\mathbb{S}^1\to \C$ is given 
by
\begin{align}
		F(\hat{\bx}) = -\int_{\Gamma}\re^{-\ri k_1\hat{\bx}\cdot\by}\left(\ri k_1(\hat{\bx}\cdot\bn(\by))u_1(\by)+\frac{\partial u_1}{\partial\bn}(\by)\right)\sd s(\by), \qquad \hat{\bx}\in\mathbb{S}^1.
		\label{eqn:2D_FF}
\end{align}

To derive a BIE from \rf{eqn:u1}--\rf{eqn:u2} one takes Dirichlet and Neumann traces of both equations and applies the standard jump relations for traces of $\cS_j$ and $\cD_j$ (see, e.g., \cite[p.~115]{acta}). This produces a set of four equations satisfied by the unknown boundary data
\[
\bv:=
  \left.\left(
    \begin{array}{c}
      u_1 \\
      \partial u_1/\partial \bn
    \end{array}
  \right)\right|_{\Gamma}
  =\left.\left(
    \begin{array}{c}
      u_2 \\
\alpha\, \partial u_2/\partial \bn
    \end{array}
  \right)\right|_{\Gamma},
\]
which can be combined in a number of different ways to produce different BIE formulations. The particular formulation we consider involves the BIE
\begin{equation}
    A\bv  = \mathbf{f},
    \label{matrix}
\end{equation}
where
\begin{align}
\label{AfDef}
 A =
 \begin{pmatrix}
  \frac{1+\alpha}{2}I+(\alpha D_2-D_1) & S_1-S_2 \\
  \alpha(H_2-H_1) &  \frac{1+\alpha}{2}I+(\alpha D_1'-D_2')
 \end{pmatrix}, \quad
 \qquad  \mathbf{f}=\left(
    \begin{array}{c}
      u^i \\
\alpha\,  \partial u^i/\partial \bn
    \end{array}
  \right),
\end{align}
$I$ is the identity operator and $S_j$, $D_j$, $D_j'$, $H_j$, $j=1,2$, are, respectively, the single-layer, double-layer, adjoint double-layer, and hypersingular integral operators, defined for sufficiently smooth $\phi$
by
\begin{eqnarray*}
   &&S_j\phi(\bx):=\int_{\Gamma}\Phi_j(\bx,\by)\phi(\by)\rd s(\by),
  \qquad \, D_j\phi(\bx):=\int_{\Gamma}\frac{\partial\Phi_j(\bx,\by)}{\partial\bn(\by)}\phi(\by)\rd s(\by), \\
   &&D'_j\phi(\bx):=\int_{\Gamma}\frac{\partial\Phi_j(\bx,\by)}{\partial\bn(\bx)}\phi(\by)\rd s(\by),\quad
H_j\phi(\bx) :=\frac{\partial}{\partial\bn(\bx)}\int_{\Gamma}\frac{\partial\Phi_j(\bx,\by)}{\partial\bn(\by)}\phi(\by)\rd s(\by).
\end{eqnarray*}
The operator $A$ in \rf{AfDef} is bounded and invertible as a mapping from $L^2(\Gamma)\times L^2(\Gamma)\to L^2(\Gamma)\times L^2(\Gamma)$; this can be proved by a straightforward modification of the arguments presented in \cite[Proof of Theorem~7.2]{torres1993helmholtz} (for details see \cite[\S2.7]{GrothPhDThesis}). 
That the solution of the BVP is sufficiently smooth to work in this setting follows from the boundedness and invertibility of $A$ as a mapping from $H^{1}(\Gamma)\times L^2(\Gamma)\to H^{1}(\Gamma)\times L^2(\Gamma)$ (again, see \cite{torres1993helmholtz,GrothPhDThesis}) and the fact that $u^i$ is smooth, so that $\mathbf{f}\in H^1(\Gamma) \times L^2(\Gamma)$.
In fact, by standard elliptic regularity results (see, e.g., \cite{GilbargTrudinger}) the solution $\bv$ is infinitely smooth on each of the sides of the polygon. Singular behaviour of $\bv$ at the corners of the polygon will be captured in our BEM using mesh refinement. 

Our choice of the particular formulation \rf{matrix}, which is similar to that used in \cite[\S3.8]{Colton} (where only smooth scatterers are considered) and also to that used in \cite{torres1993helmholtz} (albeit for an indirect method, in which the unknowns are non-physical ``densities'', rather than the boundary data itself), was made because the cancellation of the strong singularities between the two hypersingular operators in the term $H_2-H_1$ makes implementation particularly simple. 
However, we emphasise that our HNA approximation space (described in the next section) can be applied in the context of \emph{any} direct BIE formulation in which the BIE solution is some combination of the physical unknowns $u_1$ and $\pdonetext{u_1}{\bn}$, see, e.g., \cite{costabel,kleinman1988single,LaRaSa:09,Hsiao:11,DoLyTu:16}.  

We end this section with a further comment on the choice of BIE formulation for high frequency transmission problems. When using conventional (i.e.,\ non-HNA) approximation spaces, it is well-known that conditioning issues can severely limit the performance of iterative methods that are necessary for solving the large linear systems that arise. Because of this, the development of new well-conditioned BIE formulations is a highly active area of current research; in particular we mention the indirect formulations recently proposed and analysed in \cite{boubendir2015integral,boubendir2015regularized,boubendir2015high}, which offer favourable spectral properties compared to other existing formulations. (The recent paper \cite{DoLyTu:16} describes a direct counterpart with similar properties.) By contrast, we emphasise that the systems arising from our HNA approach are small enough that they can be solved \emph{directly} rather than iteratively. Indeed, the whole aim of HNA methods is to obtain small linear systems whose size does not grow with respect to frequency, hence reducing the importance of conditioning. 

\section{HNA approximation space}
\label{sec:HNA}
Our BEM solves the integral equation \rf{matrix} by the Galerkin method, using an HNA approximation space specially designed to capture the high frequency behaviour of the solution $\bv=(u_1,\pdonetext{u_1}{\bn})^T$. 
Specifically, our method is based on the ansatz
\begin{align}
\label{eqn:vDecomp}
\bv(\bx)\approx \bv^{GO}(\bx) + \bv^{d}(\bx), \qquad \bx\in\Gamma,  
\end{align}
where 
$\bv^{GO} = (u^{GO},\pdonetext{u^{GO}}{\bn})$ is the GO approximation (computed using the beam tracing algorithm described in \S\ref{sec:comp_v_go}), and $\bv^{d} = (u^{d},\pdonetext{u^{d}}{\bn})$ is a contribution from diffracted fields (computed using the HNA approximation strategy described in \S\ref{sec:Vd}). 
Before spelling out in more detail the practicalities of how $\bv^{GO}$ and $\bv^{d}$ are computed, we pause to briefly explain the origin of \rf{eqn:vDecomp} and discuss the accuracy of approximation one might expect it to provide. 

As was mentioned in the introduction, the high frequency asymptotic theory for transmission problems is not fully understood, because the details of the GTD for penetrable scatterers, in particular the solution of the canonical penetrable wedge problem, are yet to be fully worked out. 
However, by analogy with the impenetrable case, one expects the leading order high frequency behaviour of $\bv$ to be given by the GO approximation $\bv^{GO}$, which consists of the incident wave, along with reflected and transmitted waves satisfying Snell's law and the Fresnel formulas. For penetrable polygons irradiated by a plane wave, neighbouring parallel rays incident on the same side of the polygon remain parallel under reflection/transmission, and hence $\bv^{GO}$ can be computed efficiently using a beam tracing algorithm, as we describe in \S\ref{sec:comp_v_go}. As mentioned in \S\ref{sec:intro}, the field approximation obtained by substituting this GO approximation into the integral representation formulas 
\rf{eqn:u1} and \rf{eqn:u2} is referred to as the PGOH method in the atmospheric physics literature, and can be interpreted as a generalization of the classical Kirchhoff approximation for impenetrable scatterers. 

As the next-order correction to the GO approximation, we expect diffracted wave fields emanating radially outwards from each of the corners of the polygon into the exterior domain $\Omega_1$ (with wavenumber $k_1$) and inwards into the interior domain $\Omega_2$ (with wavenumber $k_2$). It is these diffracted fields, denoted $\bv^d$, that we compute numerically using our BEM, using the HNA approximation strategy described in \S\ref{sec:Vd}; see in particular the ansatz in equation \rf{eqn:ansatztransmission}.

\begin{figure}[p!]
\centering
\subfigure[Total field (real part)]{%
      \includegraphics[width=0.55\textwidth]{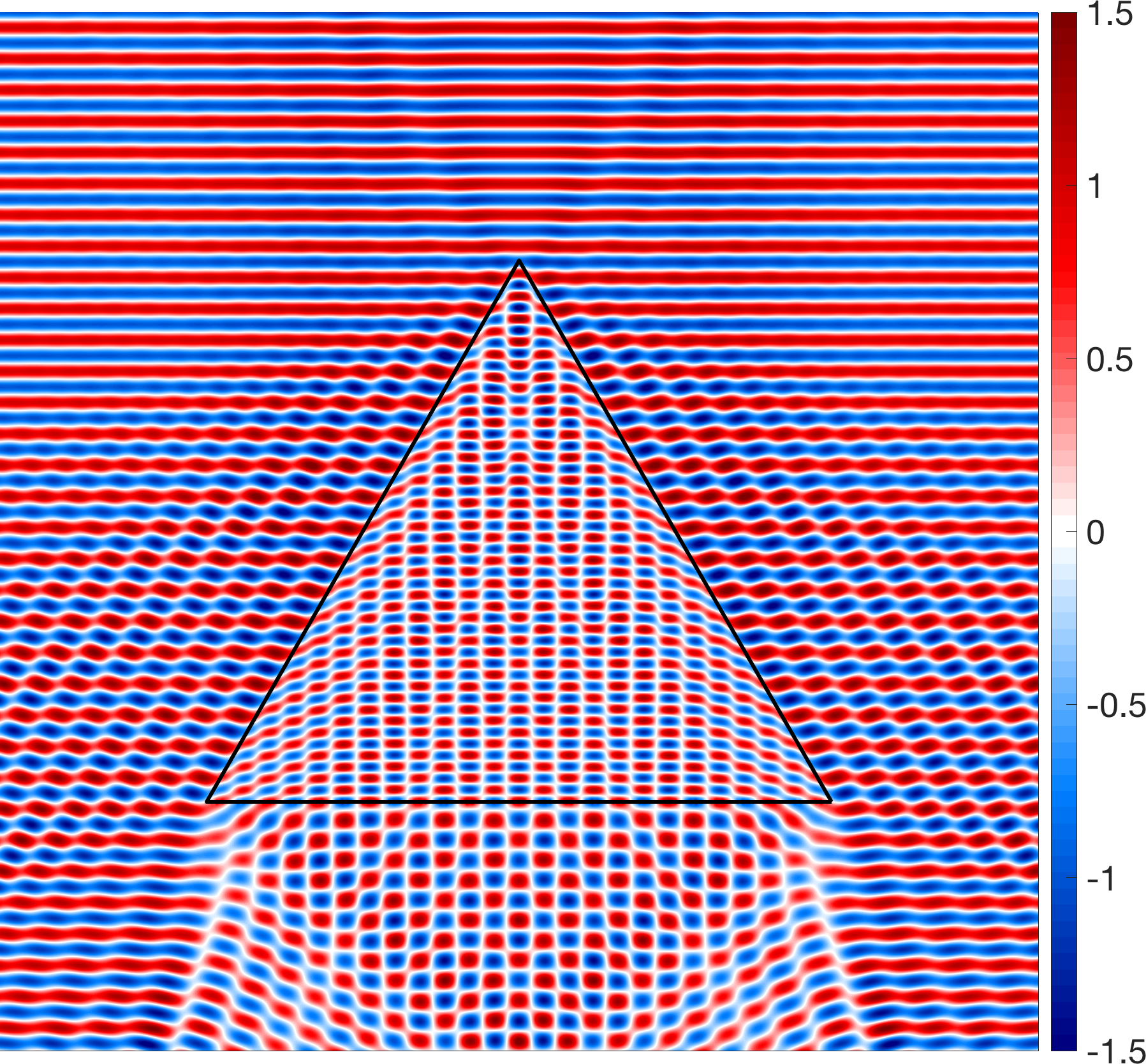}
    }    \vs{-2}
        
    \subfigure[Diffracted field (real part)]{%
      \includegraphics[width=0.55\textwidth]{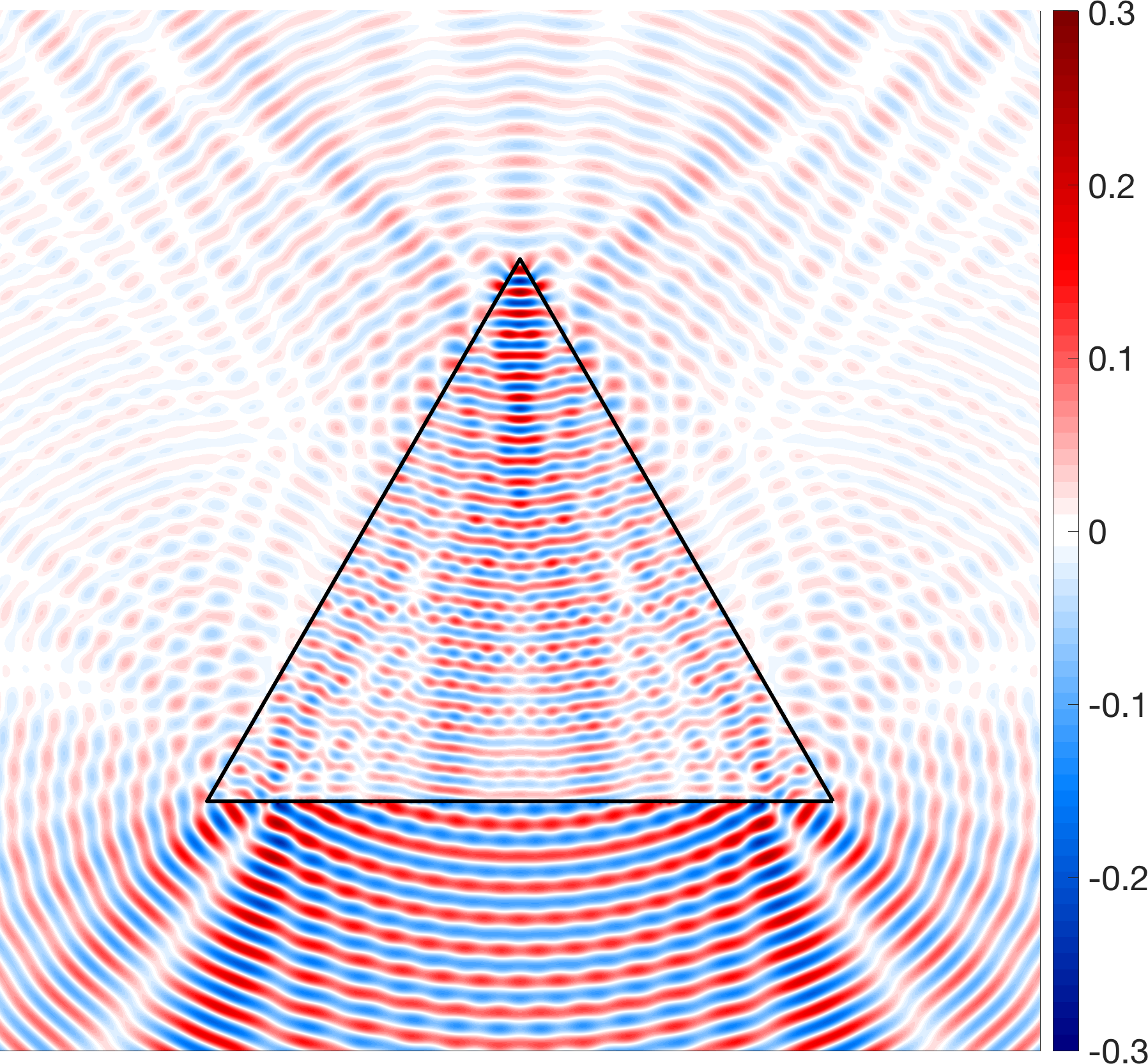}
    }
\vs{-3}
\caption{Scattering by an equilateral triangle of side length $2\pi$ and refractive index $\mu=1.5+0.003125\ri$, computed using our HNA BEM. Here $k_1=20$, $\alpha=1$, 
and the 
incident wave arrives from above, corresponding to 
$\bd^i=\bd_1$
in the geometrical setup described in Figure \ref{fig:triangle_setup}(a). 
} 
   \label{fig:tri_decomp}
\end{figure}

\begin{figure}[t!]
\centering
      \includegraphics[width=0.9\textwidth]{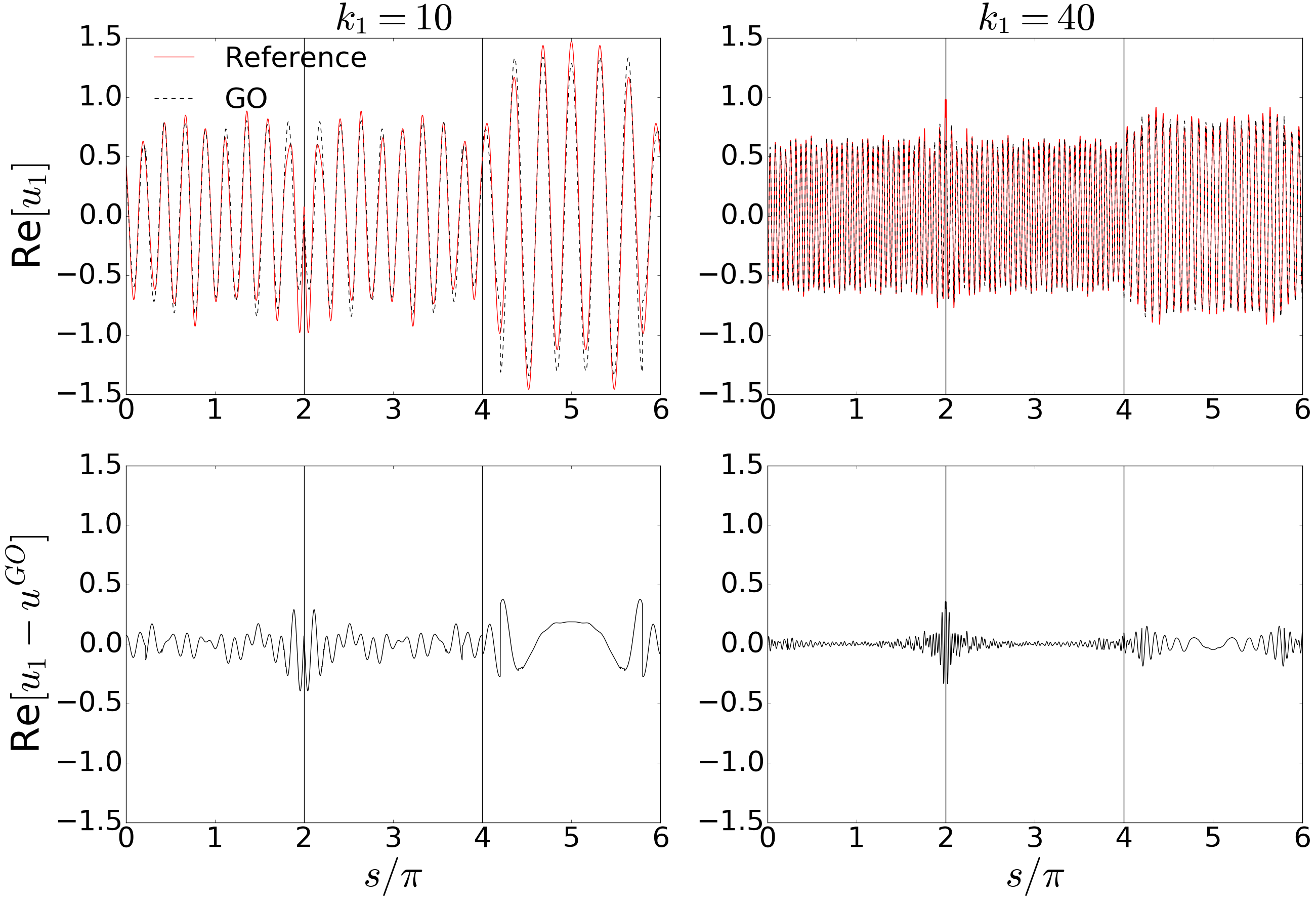}
      \caption{Real part of the boundary solution $u_1$ for the configuration in Figure~\ref{fig:tri_decomp}, with $k_1=10$ (left-hand panels) and $k_1=40$ (right-hand panels), plotted against arclength $s$ around $\Gamma$ (moving anti-clockwise around $\Gamma$ starting at the bottom-right vertex). The upper panels show the reference solution $u_1$ and the GO approximation $u^{GO}$. The lower panels show the difference $u_1-u^{GO}$ which we identify with the diffracted field. This difference, rather than the total field, is what we approximate using our HNA BEM. The vertical lines indicate the location of the corners of the polygon. 
      }
      \label{fig:boundary_decomp}
\end{figure}

The contribution made by the diffracted term $\bv^d$ is illustrated in Figures~\ref{fig:tri_decomp} and \ref{fig:boundary_decomp}, which relate to scattering by an equilateral triangle of refractive index $\mu=1.5+0.003125\ri$ (for full details of the setup see \S\ref{sec:Numerical}). Figure~\ref{fig:tri_decomp}(a) shows the real part of the total field computed using our HNA BEM, i.e., first calculating the approximation $\bv\approx \bv^{GO}+\bv^d$ and then substituting this into the representation formulas \rf{eqn:u1} and \rf{eqn:u2}. Figure~\ref{fig:tri_decomp}(b) shows the contribution made by the diffracted term $\bv^d$ alone; here one can clearly see the circular diffracted waves emanating from the vertices, and the complicated interference pattern they produce inside the scatterer.  
\dhednote{I wasn't sure what else to add here. To be honest, I would be happy to delete Figure~\ref{fig:boundary_decomp} completely, to speed up the flow, and avoid having to explain it. But I know Sam is fond of it...}
\sgednote{Up to you. I like this picture since I think it shows the $u_{go}+u_d$ decomposition. I've added a sentence to the caption to highlight that the difference $u-u_{go}$ is what we actually approximate. This is something which I feel Reviewer 1 didn't understand.}
Figure~\ref{fig:boundary_decomp} contains plots of the solution $u_1$ on $\Gamma$, i.e.,\ the first component of the vector $\bv$. 
The upper panels show the real part of $u_1$ on $\Gamma$, computed both using a highly accurate reference solution (described in \S\ref{sec:Numerical}) and using the GO approximation $u^{GO}$. 
The lower panels show the difference $u_1-u^{GO}$; this remainder is what we approximate using our diffracted term $u^d$. 
A comparison of the left-hand panels ($k_1=10$) with the right-hand panels ($k_1=40$), suggests that as the wavenumber $k_1$ increases the error $u_1-u^{GO}$ decreases in magnitude, 
at least away from the corners of the polygon. (Similar behaviour is observed for the GO approximation to $\pdonetext{u_1}{\bn}$, not shown here.) This is in line with what one should expect from the high frequency GO approximation, i.e., that it becomes more accurate with increasing $k_1$, but that it does not capture any effects due to the corner diffraction.

\begin{figure}[t!]
\centering
\subfigure[Lateral waves]{\includegraphics[width=65mm]{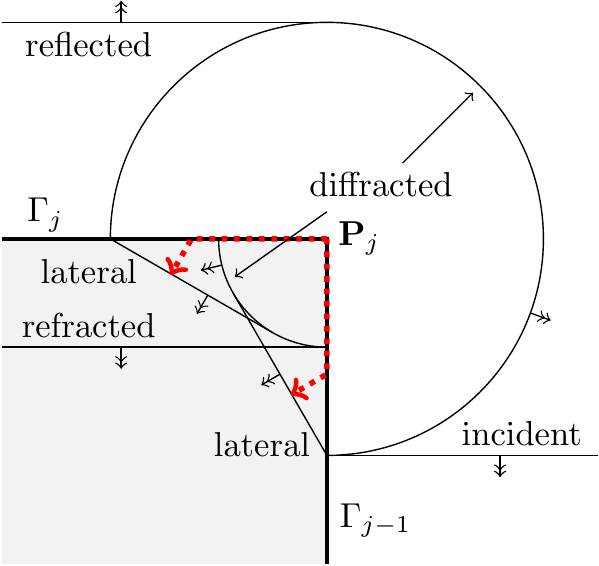}
}
\hs{20}
\subfigure[Diffracted-reflected wave]{\includegraphics[width=55mm]{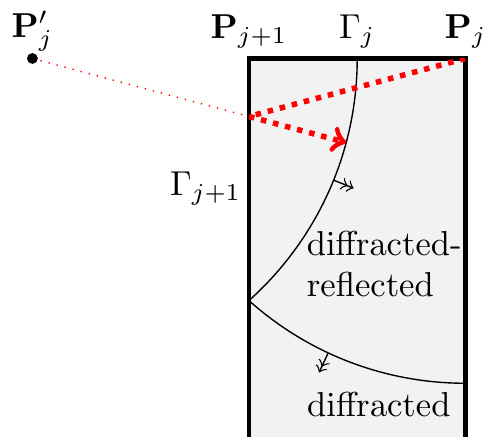}
}
\caption{Wavefront propagation in the time domain for a plane pulse incident from a faster medium onto a rectangle formed of a slower medium (shaded), illustrating two of the high frequency asymptotic phenomena \emph{not} captured by our HNA approximation space. Double-headed arrows indicate wavefront propagation directions. (a) Wavefront diagram near a corner of the rectangle shortly after the incident pulse (which propagates downwards in the figure) arrives at the corner $\bP_j$, showing the lateral wavefronts, which are associated with the limiting exterior diffracted rays propagating down the sides $\Gamma_{j-1}$ and $\Gamma_{j}$ shedding new lateral rays into the interior. Typical lateral ray paths are shown as dotted red arrows. In this particular ``grazing incidence'' configuration the limiting exterior diffracted ray propagating down the side $\Gamma_{j-1}$ coincides with the limiting incident ray. As a result we might expect a stronger lateral wave contribution from $\Gamma_{j-1}$ than in the general case, which is illustrated, for example in \cite[Figure 3]{groth2015hybrid}. We note that the two lateral wavefronts meet the interior (slower) diffracted wavefront tangentially. (b) Diffracted wavefront emanating from the corner $\bP_j$ and the diffracted-reflected wavefront generated by its interior re-reflection from the side $\Gamma_{j+1}$. A typical diffracted-reflected ray path is shown by the dotted red arrow.
\label{fig:LateralAndDiffRef}}
\end{figure}

We also expect the high frequency solution behaviour to feature other higher-order effects, such as so-called \emph{lateral }waves (sometimes known as \emph{head} or \emph{bow} waves in related contexts), and the internal re-reflections of the diffracted and lateral wave fields. (For a more detailed discussion see \cite[\S3.2.1]{groth2015hybrid}, and for a schematic illustration of some of these phenomena see Figure \ref{fig:LateralAndDiffRef}.) 
However, to limit the implementational complexity of our method these higher-order effects are not incorporated into our approximation space. 
Hence, in contrast to the HNA approximation spaces developed in \cite{CWL:07,CWLM,hewett2013high} for impenetrable convex polygons, our HNA approximation space does not capture \emph{all} of the oscillatory solution behaviour. Accordingly, our method remains ``asymptotic'' in nature, with an inherent frequency-dependent error which decreases as the frequency tends to infinity. \footnote{Many HNA methods share this same ``asymptotic'' nature. In particular we note that the methods of \cite{dominguez2007hybrid,bruno2004prescribed,GH11} for smooth convex impenetrable scatterers all incur an exponentially small frequency-dependent error, since they approximate by zero in the deep shadow region, neglecting the exponentially small creeping wave fields. 
In contrast, the method of \cite{Gi:07} and the analysis of \cite{AsHu:14} do explicitly incorporate creeping wave fields, and the methods of \cite{EcOz:17,EcUr:16}, while not fully capturing creeping wave fields, include enough degrees of freedom in the deep shadow to be fully convergent, without incurring a systematic frequency-dependent error.}
Nonetheless, our numerical results in \S\ref{sec:Numerical} show that our method achieves a significant improvement over using GO alone. 
As we shall see in \S\ref{sec:Numerical}, our neglect of the higher-order terms (lateral waves, and the re-reflections of the diffracted and lateral waves) has less impact when the scatterer is partially absorbing ($\im{k_2}>0$), since then the neglected terms are attenuated and play less of a role. 

We now describe the computation of each of the two terms $\bv^{GO}$ and $\bv^d$ in more detail.

\subsection{Computation of $\bv^{GO}$}
\label{sec:comp_v_go}
Our GO approximation $\bv^{GO}$ in \rf{eqn:vDecomp} is computed 
using a beam-tracing algorithm (BTA), similar to those presented in \cite{borovoi2003scattering} and \cite{Bi:scat} for the analogous 3D problem. For full details of our BTA we refer the reader to \cite[\S3.1]{groth2015hybrid}; 
here we simply sketch the basic algorithm and point out some modifications in our current implementation compared to that described in \cite{groth2015hybrid}.

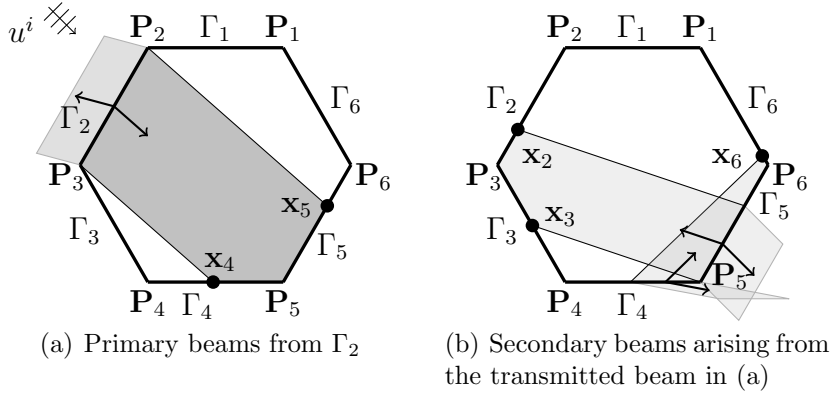
\begin{figure}[t!]
\centering
\subfigure[Primary beams from $\Gamma_2$]
{
\begin{tikzpicture}[scale=0.6,baseline]
\tikzstyle{conefill} = [fill=black!30,fill opacity=0.8]
\tikzstyle{ann} = [fill=white,font=\footnotesize,inner sep=1pt]   
    \filldraw[fill=black!30,draw=black!30!,fill opacity=0.4]
    (-1.5,2.5981)--(-3,0)--(-3.9659,0.2588)
              --(-2.4659,2.8569)--cycle;  
                       
    \draw[->,black,thick](-2.25,1.299) -- (-3.1193,1.5319);     
   \tikzstyle{conefill} = [fill=black!40,fill opacity=0.6]     
   \filldraw[conefill](-1.5,2.5981)--(-3,0)
                       --(-0.0525,-2.5981)--(1.5,-2.5981)--(2.4764,-0.9070)--cycle;   
    \draw[->,thick](-2.25,1.299) -- (-1.4998,0.6378);  
    \draw[line width=1.25pt](1.5,2.5981)--(-1.5,2.5981) ;
    \draw[line width=1.25pt](-1.5,2.5981)--(-3,0) ;
    \draw[line width=1.25pt](-3,0)--(-1.5,-2.5981) ;
    \draw[line width=1.25pt](-1.5,-2.5981)--(1.5,-2.5981) ;
    \draw[line width=1.25pt](1.5,-2.5981)--(3,0) ;
    \draw[line width=1.25pt](3,0)--(1.5,2.5981) ;
    \draw(0,2.5) node[anchor=south]{$\Gamma_1$};
    \draw(-0.4,-2.5981) node[anchor=north]{$\Gamma_4$};
    \draw(-2.5,1.0) node[anchor=east]{$\Gamma_2$};
    \draw(-2.3,-1.45) node[anchor=east]{$\Gamma_3$};
    \draw(2.3,1.45) node[anchor=west]{$\Gamma_6$};
    \draw(2.0,-1.8) node[anchor=west]{$\Gamma_5$};
    \draw(1.5,2.5) node[anchor=south]{$\bP_1$};
        \draw(-1.5,2.5) node[anchor=south]{$\bP_2$};    
                \draw(-3.3,-0.8) node[anchor=south]{$\bP_3$};    
                \draw(-1.5,-3.6) node[anchor=south]{$\bP_4$};
            \draw(1.5,-3.6) node[anchor=south]{$\bP_5$};
	\draw(3.5,-0.8) node[anchor=south]{$\bP_6$};    
	\draw(0.1,-1.7) node[anchor=north]{$\bx_4$};
	\draw(1.8,-0.5) node[anchor=north]{$\bx_5$};
     \draw[->]
          (-3.7,3.5)--(-3.1,2.9);
          \draw (-4.3,3) node{$u^i$};      
     \draw[black](-3.8,3.2)--(-3.4,3.6);
    \draw[black](-3.5,2.9)--(-3.1,3.3);
    \draw[black](-3.65,3.05)--(-3.25,3.45);    
        \fill[black] (2.4764,-0.9070) circle (0.15); 
    \fill[black] (-0.0525,-2.5981) circle (0.15); 
\end{tikzpicture}
}
\hspace{0.15cm}
\subfigure[Secondary beams arising from the transmitted beam in (a)]
{
\begin{tikzpicture}[scale=0.6,baseline]
\tikzstyle{conefill} = [fill=black!30,fill opacity=0.8]
\tikzstyle{ann} = [fill=white,font=\footnotesize,inner sep=1pt]       
    \tikzstyle{conefill} = [fill=black!30,fill opacity=0.2]    
    \filldraw[conefill](-0.0525,-2.5981)--(1.5,-2.5981)--(3,0)
                       --(2.87,0.2)--cycle;    	
    
    \tikzstyle{conefill} = [fill=black!30,draw=black!30,fill opacity=0.2]                  
    \filldraw[conefill](-0.0525,-2.5981)--(1.5,-2.5981)--(3.4654,-2.9683)
                       --(1.9129,-2.9683)--cycle;

   \tikzstyle{conefill} = [fill=black!30,fill opacity=0.2]
        
   \filldraw[conefill](-2.5473,0.7842)--(-3,0)--(-2.2237,-1.3445)
                       --(1.5,-2.5981)--(2.4764,-0.9070)--cycle;

   \tikzstyle{conefill} = [fill=black!30,draw=black!30,fill opacity=0.2]    
   \filldraw[conefill](3.3249,-1.7555)--(2.3485,-3.4466)
                       --(1.5,-2.5981)--(2.4764,-0.9070)--cycle; 
   \draw[->,thick](0.7237,-2.5981) -- (1.4,-1.9368);                   
                       
   \draw[->,thick](1.9882,-1.7525) -- (1.0404,-1.4335);    
   
   \draw[->,thick,black](0.7237,-2.5981)--(1.7064,-2.7832);
   
   \draw[->,thick,black](1.9882,-1.7525) -- (2.6953,-2.4596); 
    \draw[line width=1.25pt](1.5,2.5981)--(-1.5,2.5981) ;
    \draw[line width=1.25pt](-1.5,2.5981)--(-3,0) ;
    \draw[line width=1.25pt](-3,0)--(-1.5,-2.5981) ;
    \draw[line width=1.25pt](-1.5,-2.5981)--(1.5,-2.5981) ;
    \draw[line width=1.25pt](1.5,-2.5981)--(3,0) ;
    \draw[line width=1.25pt](3,0)--(1.5,2.5981) ;
    \draw(0,2.5) node[anchor=south]{$\Gamma_1$};
    \draw(0,-2.5981) node[anchor=north]{$\Gamma_4$};
    \draw(-2.3,1.45) node[anchor=east]{$\Gamma_2$};
    \draw(-2.3,-1.45) node[anchor=east]{$\Gamma_3$};
    \draw(2.3,1.45) node[anchor=west]{$\Gamma_6$};
    \draw(2.55,-0.9) node[anchor=west]{$\Gamma_5$};    
        \draw(1.5,2.5) node[anchor=south]{$\bP_1$};
        \draw(-1.5,2.5) node[anchor=south]{$\bP_2$};    
                \draw(-3.3,-0.8) node[anchor=south]{$\bP_3$};    
                \draw(-1.5,-3.6) node[anchor=south]{$\bP_4$};
            \draw(2.15,-3) node[anchor=south]{$\bP_5$};
	\draw(3.5,-0.8) node[anchor=south]{$\bP_6$};    
	\draw(-1.6,-0.7) node[anchor=north]{$\bx_3$};
	\draw(-2.1,0.6) node[anchor=north]{$\bx_2$};
		\draw(2.1,0.6) node[anchor=north]{$\bx_6$};
         \fill[black] (2.87,0.2) circle (0.15); 
                                       \fill[black] (-2.5473,0.7842) circle (0.15); 
                              \fill[black] (-2.2237,-1.3445) circle (0.15); 
\end{tikzpicture}
}
\caption{Beam tracing to compute $\bv^{GO}$ in a hexagon: (a) shows the primary reflected and transmitted beams arising from the incidence of $u^i$ onto the side $\Gamma_2$; (b) shows the secondary beams arising from the internal reflection and exterior transmission of the primary transmitted beam in (a). The secondary internally-reflected beams in (b) go on to produce a sequence of higher-order internally-reflected beams, not pictured here. For the particular incident direction illustrated here, the sides $\Gamma_1$ and $\Gamma_3$ would also produce a similar sequence of beams (not pictured here) contributing to $\bv^{GO}$, but $\Gamma_4$, $\Gamma_5$ and $\Gamma_6$ would not, as they are not illuminated by $u^i$. 
The beam boundary points associated with the beams shown here are indicated by the filled circles and labelled $\bx_j$, $j=2,\ldots,6$ to indicate the side of $\Gamma$ on which they lie; note that the secondary beam leaving $\Gamma_4$ in (b) has only one beam boundary inside the polygon, so generates only one beam boundary point, $\bx_6$.}
\label{fig:rays}
\end{figure}

According to the principles of GO, when the incident plane wave impinges on a side of the polygon $\Omega_2$, it generates a beam of reflected rays propagating back into the exterior domain $\Omega_1$, and a beam of transmitted rays propagating into the interior domain $\Omega_2$. The transmitted beam then generates a sequence of higher-order internally reflected beams and externally transmitted beams, as illustrated in Figure \ref{fig:rays}. Each beam is described by 
\begin{enumerate}[(i)]
\item \label{i} a plane wave $a\re^{\ri (\NN\dr + \ri \KK\di)\cdot\bx}$, with amplitude $a$, unit propagation and decay direction vectors $\dr$ and $\di$, and (wavenumber-dependent) propagation/decay constants $\NN>0$ and $\KK\geq 0$, all determined by the well-known laws of reflection and refraction at an interface between two homogeneous media (Snell's law \rf{eqn:Snell} and the Fresnel formulas \rf{eqn:Fresnel1}--\rf{eqn:Fresnel2} below), and 
\item \label{ii} (up to) two beam boundaries outside of which we cut off the plane wave beam sharply to zero.
\end{enumerate}

In principle there is an infinite sequence of internal reflections to consider, but for a practical implementation we need to truncate the sequence somehow. Our truncation criterion is as follows. Suppose that a beam $a\re^{\ri (\NN\dr + \ri \KK\di)\cdot\bx}$ (with associated beam boundaries) illuminates a portion of a side $\Gamma_j$. Then we include its contribution to $\bv^{GO}$ on $\Gamma_j$ if and only if 
\begin{align}
\label{tolBDef}
\max_{\bx}|a|\re^{-\KK\di\cdot\bx}> {\rm tol}_{b},
\end{align}
where the maximum is taken over all $\bx\in\Gamma_j$ which lie between the two beam boundaries of the beam in question, and ${\rm tol}_{b}$ is a user-specified tolerance. If a beam fails this test, we exclude it from the calculation of $\bv^{GO}$ on $\Gamma_j$ and do not track any higher-order beams generated by the reflection of that beam from the side $\Gamma_j$. 
The exact number of beams traced depends upon ${\rm tol}_{b}$, frequency, absorption, geometry, and incident angle. Generally speaking, higher absorption typically leads to fewer beams, whereas a more complex geometry (more sides) leads to more beams. 
For the numerical results in \S\ref{sec:Numerical} we used ${\rm tol}_{b} = 0.005$; this was found to guarantee convergence of the BTA to within a relative accuracy of around $0.1\%$ for all the examples considered, using between $5$ and $30$ beams, depending on the example. 
We emphasize that tracing each beam requires us to trace just two geometric rays - the two rays forming the beam boundaries. As a result, in all of the examples we consider, the computational cost of the BTA (both in terms of memory and computation time) is negligible within the overall cost of the HNA BEM, even when a relatively large number of beams is required (such as for the hexagon in \S\ref{subsec:geom}).

\subsubsection{Sign choice in the Fresnel formulas}
\label{sss:sign}

Regarding point \rf{i} above, we remark that the laws of reflection and refraction are completely classical in the case when both media are non-absorbing (see, e.g.,~\cite{Born}). However, when one of the media is absorbing (here, this corresponds to $\im{k_2}>0$), there is some debate in the literature about how to make a particular sign choice in the Fresnel formulas. We reviewed this issue in \cite[Appendix A]{groth2015hybrid}, providing an analysis of the general interface problem (with up to two absorbing media) as well as numerous references to relevant literature. But since the publication of \cite{groth2015hybrid}, further investigations have led us to refine our sign choice rule, as we now explain.

Let us consider the interface problem between a medium with wavenumber $k=\eta+\ri \xi$ and a second medium with wavenumber $\tilde k=\tilde\eta+\ri \tilde\xi$, where $\eta,\tilde\eta>0$, $\xi,\tilde\xi\geq0$. 
Let $\bt$ and $\bn$ be unit tangent and normal vectors on the interface (assumed to contain the origin), with $\bn$ pointing into the second medium. 
An incident plane wave \mbox{$u^i=a^i \exp[\ri (D_i \bd^i + \ri E_i \be^i)\cdot \bx]$} propagating in the first medium generates a reflected wave \mbox{$u^r=a^r \exp[\ri (D_i \bd^r + \ri E_i \be^r)\cdot \bx]$} (also propagating in the first medium), and a transmitted wave 
\mbox{$u^t=a^t \exp[\ri (D_t \bd^t + \ri E_t \be^t)\cdot \bx]$} (propagating in the second medium).
According to the law of reflection the vectors $\bd^r$ and $\be^r$ should satisfy 
\[ 
\bd^r = \bd^i - 2(\bd^i\cdot \bn)\bn, \qquad \be^r = \be^i - 2(\be^i\cdot \bn)\bn,
\]
and for $u^i$, $u^r$ and $u^t$ to be solutions of the relevant Helmholtz equations we require
\begin{align}
\label{InterfaceRelns}
D_i^2-E_i^2=\eta^2-\xi^2, 
\quad
D_t^2-E_t^2=\tilde\eta^2-\tilde\xi^2, 
\quad
D_iE_i(\bd^i\cdot \be^i)= 
\eta\xi, \quad
D_tE_t(\bd^t\cdot \be^t)=\tilde\eta\tilde\xi.
\end{align}

Assuming the boundary conditions $ u^i + u^r = u^t $ and $\pdonetext{u^i}{\bn}+\pdonetext{u^r}{\bn}=\alpha \pdonetext{u^t}{\bn}$ on the interface, for some $\alpha\in \C$, one can derive Snell's law
\begin{align}
\label{eqn:Snell}
D_i\bd^i\cdot\bt = D_t\bd^t\cdot\bt,\qquad E_i\be^i\cdot\bt = E_t\be^t\cdot\bt,
\end{align}
and the Fresnel formulas
\begin{align}
\label{eqn:Fresnel1}
	\frac{a^r}{a^i} &= \frac{(D_i\bd^i\cdot\bn+\ri E_i\be^i\cdot\bn) - \alpha(D_t\bd^t\cdot\bn + \ri E_t\be^t\cdot\bn)}{(D_i\bd^i\cdot\bn+\ri E_i\be^i\cdot\bn) + \alpha(D_t\bd^t\cdot\bn + \ri E_t\be^t\cdot\bn)}, \\
	\frac{a^t}{a^i} &= \frac{2(D_i\bd^i\cdot\bn+\ri E_i\be^i\cdot\bn)}{(D_i\bd^i\cdot\bn+\ri E_i\be^i\cdot\bn) + \alpha(D_t\bd^t\cdot\bn + \ri E_t\be^t\cdot\bn)},
\label{eqn:Fresnel2}
\end{align}
in which
\begin{align*} 
D_t &= \sqrt{\frac{1}{2}\left(\tilde{\eta}^2-\tilde{\xi}^2+\Ns^2+\Ks^2 +\sqrt{(\tilde{\eta}^2-\tilde{\xi}^2 -\Ns^2+\Ks^2)^2 +4(\Ns\Ks-\tilde{\eta}\tilde{\xi})^2}\right)},\\
         \KK _t &= \sqrt{\NN_t^2+\tilde{\xi}^2-\tilde{\eta}^2},  
\end{align*}
with $\Ns:=\NN_i\dri\cdot\bt$ and $\Ks:=\KK _i\dii\cdot\bt$ (see the appendix of \cite{groth2015hybrid} for more detail). Snell's law and the Fresnel formulas respectively specify the tangential components of the transmitted vectors $\bd^t$ and $\be^t$, and the amplitudes $a^r$ and $a^t$. Since $\bd^t$ and $\be^t$ have unit length, their normal components are thus specified up to sign. However, to determine these signs and close the problem one must impose an additional constraint based on some physical argument. 
Two common constraints require either that
\begin{align}
\label{GO1}
\bd^t\cdot \bn \geq0, \qquad \textrm{(GO1)},
\end{align}
which prohibits energy flow from the second medium to the first,
or that 
\begin{align}
\label{GO2}
\be^t\cdot \bn \geq0, \qquad \textrm{(GO2)},
\end{align}
which prohibits exponential growth in the second medium. 

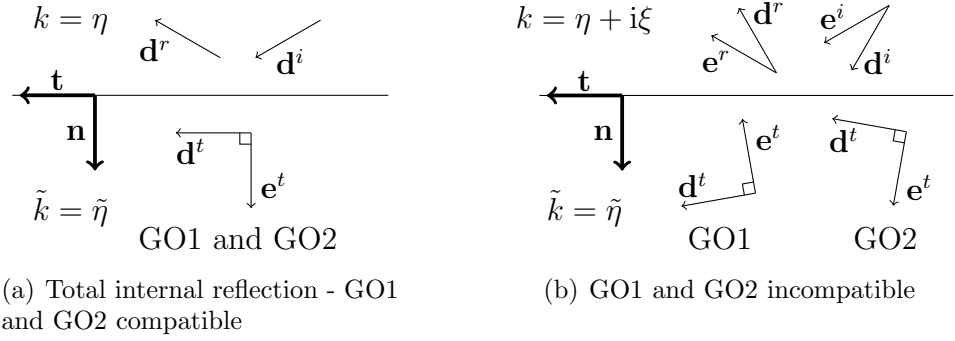
\begin{figure}[t]
\centering
\subfigure[Total internal reflection - GO1 and GO2 compatible]
{
\begin{tikzpicture}
\draw (-3,0)--(2,0);
\draw[<-,line width=0.5mm] (-2.9,0)--(-1.9,0);
\draw[<-,line width=0.5mm] (-1.9,-1)--(-1.9,0);
\draw (-2.4,0.2) node{$\bt$};
\draw (-2.15,-0.5) node{$\bn$};
\draw[->] (1.1,1)--({1.1-cos(30)},{1-sin(30)});
\draw[<-] (-1.1,1)--({-1.1+cos(30)},{1-sin(30)});
\draw (0.7,0.45) node{$\bd^i$};
\draw (-1.1,0.6) node{$\bd^r$};
\begin{scope}[xshift=5]
\draw[->] (0,-0.5)--(-1,-0.5);
\draw[->] (0,-0.5)--(0,-1.5);
\draw (-0.15,-0.5)--(-0.15,-0.65)--(0,-0.65);
\draw (-0.8,-0.75) node{$\bd^t$};
\draw (0.3,-1.2) node{$\be^t$};
\draw (-2.4,1.0) node{$k=\eta$};
\draw (-2.4,-1.5) node{$\tilde k=\tilde{\eta}$};
\end{scope}
\draw (0,-1.9) node{GO1 and GO2};
\end{tikzpicture}
\label{TIR}
}
\hs{10}
\subfigure[GO1 and GO2 incompatible]
{
\begin{tikzpicture}
\draw (-3,0)--(2.5,0);
\draw[<-,line width=0.5mm] (-2.9,0)--(-1.9,0);
\draw[<-,line width=0.5mm] (-1.9,-1)--(-1.9,0);
\draw (-2.4,0.2) node{$\bt$};
\draw (-2.15,-0.5) node{$\bn$};
\begin{scope}[xshift=10]
\draw[->] (1.3,1.2)--({1.3-cos(30)},{1.2-sin(30)});
\draw[->] (1.3,1.2)--({1.3-cos(60)},{1.2-sin(60)});
\draw[->] (-0.2,0.3)--({-0.2-cos(60)},{0.3+sin(60)});
\draw[->] (-0.2,0.3)--({-0.2-cos(30)},{0.3+sin(30)});
\draw (1.15,0.45) node{$\bd^i$};
\draw (0.55,1.05) node{$\be^i$};
\draw (-0.3,1.1) node{$\bd^r$};
\draw (-1,0.45) node{$\be^r$};
\end{scope}
\begin{scope}[rotate=-10,xshift=55,yshift=10]
\draw[->] (0,-0.5)--(-1,-0.5);
\draw[->] (0,-0.5)--(0,-1.5);
\draw (-0.8,-0.75) node{$\bd^t$};
\draw (0.3,-1.2) node{$\be^t$};
\draw (-0.15,-0.5)--(-0.15,-0.65)--(0,-0.65);
\end{scope}
\begin{scope}[rotate=10,yscale=-1,xshift=-10,yshift=50]
\draw[->] (0,-0.5)--(-1,-0.5);
\draw[->] (0,-0.5)--(0,-1.5);
\draw (-0.8,-0.75) node{$\bd^t$};
\draw (0.3,-1.2) node{$\be^t$};
\draw (-0.15,-0.5)--(-0.15,-0.65)--(0,-0.65);
\end{scope}
\draw (-2.4,1.0) node{$k=\eta+\ri \xi$};
\draw (-2.4,-1.5) node{$ \tilde k=\tilde{\eta}$};
\draw (-0.6,-1.9) node{GO1};
\draw (1.6,-1.9) node{GO2};
\end{tikzpicture}
\label{GO1vsGO2}
}
\caption{Examples illustrating the sign choice for the normal components of transmitted direction vectors in the canonical GO interface problem. In both figures the interface is the horizontal line and the vectors $\bt$ and $\bn$ are unit tangent and normal vectors to the interface. (a) shows total internal reflection, which occurs when $0<\tilde{\eta}<\eta$ and the incident angle is sufficiently shallow (i.e., $|\bd^i\cdot \bt|$ is sufficiently small); (b) shows the case where the choices GO1 and GO2 are incompatible, which occurs when $\eta,\tilde{\eta}>0$, $\xi>0$ and $(\bd^i\cdot\bt)(\be^i\cdot\bt)>0$.}
\label{fig:vectors}
\end{figure}

For some configurations these two choices are compatible, and can both be satisfied simultaneously. As an example we consider the phenomenon of total internal reflection, which occurs when both media are non-absorbing ($\xi=\tilde\xi=0$), the second medium has a faster propagation speed than the first (i.e.,  $0<\tilde\eta<\eta$) and the incident direction is close to grazing (i.e., $|\bd^i\cdot\bt|$ is sufficiently small). In this case $\bd^t\cdot \bn=0$, so GO1 is automatically satisfied. But $\be^t$ can be either $\bn$ or $-\bn$. Enforcing GO2 selects the first option, and we obtain the configuration shown in Figure \ref{fig:vectors}(a). 

However, when one or both of the media are absorbing, there exist configurations for which GO1 and GO2 are incompatible. (Note that $\bd^t$ and $\be^t$ are coupled through the final equation in \rf{InterfaceRelns} so the signs of $\bd^t\cdot \bn$ and $\be^t\cdot \bn$ cannot in general be chosen independently.) 
For the case when just one of the media is absorbing (as in the current paper), the only situation in which this incompatibility arises is when the first medium is absorbing ($\xi>0$) and the second is non-absorbing ($\tilde\xi=0$), and the incident vectors $\bd^i$ and $\be^i$ both point in the same direction relative to $\bt$ (i.e., $(\bd^i\cdot\bt)(\be^i\cdot\bt)>0$). In this case Snell's law implies that $(\bd^t\cdot\bn)(\be^t\cdot\bn)<0$, so that \rf{GO1} and \rf{GO2} cannot both hold simultaneously. The two possibilities are illustrated in Figure \ref{fig:vectors}(b). 

In the context of our problem, this situation can arise when the scatterer is absorbing ($\im{k_2}>0$) and a beam propagating in $\Omega_2$ is incident on $\Gamma$, generating an internally reflected beam and an externally transmitted beam. We expect that the correct choice of GO1 or GO2 in this context could in principle be resolved by matching the local asymptotic behaviour near the interface with the global behaviour of the scattered field. However, due to the lack of a convenient expression for the field near the corners of the polygon we have been unable to carry out such an analysis. 
Instead we have developed a heuristic approach for determining which choice to make, based on the results of numerical experiments.  In these experiments (detailed more fully in \cite[\S4.3.4]{GrothPhDThesis}) the GO approximation $\bv^{GO}$, computed using first GO1 and then GO2, was compared to a conventional BEM reference solution, for a particular geometrical configuration chosen specifically to isolate the behaviour in question. We found that GO1 gave a more accurate approximation than GO2, except when $\bd^t$ was almost parallel to the interface. On the basis of these results we adopt the following strategy: whenever GO1 and GO2 are incompatible (as described above),
\begin{align}
\label{tolGODef}
\textrm{use GO1, unless } |\bd^t\cdot\bn|<{\rm tol}_{GO}, \textrm{ in which case use GO2}.
\end{align}
Here ${\rm tol}_{GO}$ is a user-defined tolerance; for the numerical results in \S\ref{sec:Numerical} we used the value ${\rm tol}_{GO}=0.01$. 
(In our previous paper \cite{groth2015hybrid} we used GO1 throughout, which corresponds to setting ${\rm tol}_{GO}=0$. The significant improvement in accuracy resulting from the apparently minor change to ${\rm tol}_{GO}=0.01$, enabling GO2 to be used where appropriate, is described in \cite[\S4.3.4]{GrothPhDThesis}.)

\subsubsection{Beam boundary discontinuities\label{sec:BeamBoundaries}}

Regarding point \rf{ii} above, we note that cutting off each beam sharply across its two beam boundaries introduces artificial discontinuities in our approximation \rf{eqn:vDecomp} to the function $\bv$ at the ``beam boundary points'' where the beam boundaries intersect $\Gamma$ (see Figure \ref{fig:rays}). If these discontinuities are sufficiently large in magnitude they can cause difficulties for our numerical calculation of the diffracted component $\bv^d$. We discuss this issue in more detail in the next section, and for now we simply introduce some terminology. Suppose that a beam $a\re^{\ri (\NN\dr + \ri \KK\di)\cdot\bx}$ (with associated beam boundaries) impinges on a side $\Gamma_j$, with one of its beam boundaries intersecting $\Gamma_j$ at the point $\bx_*\in\Gamma_j$. Then we designate the resulting beam boundary discontinuity as ``strong'' if
\begin{align}
\label{tolBBDef}
|a|\re^{-\KK\di\cdot\bx_*}> {\rm tol}_{bb},
\end{align}
where 
${\rm tol}_{bb}$ is a user-specified tolerance.

\subsection{Computation of $\bv^d$}
\label{sec:Vd}

Our diffraction contribution $\bv^d$ in \rf{eqn:vDecomp} takes the form
\begin{align}
\label{eqn:ansatztransmission}
\bv^{d}(\bx) = 
\sum_{j=1}^{n_s}
\bv_1^{j}(\bx) \exp[\ri k_1 r_{j}(\bx)]
+\sum_{j=1}^{n_s}
\bv_2^{j}(\bx)\exp[\ri k_2 r_{j}(\bx)],
\qquad \bx\in\Gamma,
\end{align}
where, for $j=1,\ldots,n_s$, $r_j(\bx)=|\bx-\bP_j|$ is the Euclidean distance between $\bx\in\Gamma$ and the corner $\bP_j$, and the amplitudes $\bv_l^j$, $l=1,2$, are computed numerically using our Galerkin BEM (see~(\ref{eqn:gal}) below), as (vector-valued) piecewise polynomials on appropriate overlapping meshes. (Recall that since $\bv=(u_1,\pdonetext{u_1}{\bn})$ is vector-valued, so are the amplitudes $\bv_l^j$.)

A typical term $\bv_1^{j}(\bx) \exp[\ri k_1 r_{j}(\bx)]$ in the first sum in \rf{eqn:ansatztransmission} corresponds to a diffracted wave emanating from the corner $\bP_j$ and propagating with wavenumber $k_1$ in the exterior domain $\Omega_1$. The convexity of the polygon $\Omega_2$ implies that this wave only impinges on the sides adjacent to the corner $\bP_j$, namely $\Gamma_j$ and $\Gamma_{j-1}$ (with the convention that $\Gamma_0\equiv \Gamma_{n_s}$). Hence we force $\bv_1^j(\bx)=\boldsymbol{0}$ for $\bx\in\Gamma\setminus(\Gamma_j\cup\Gamma_{j-1})$. 
On each of $\Gamma_j$ and $\Gamma_{j-1}$, to capture the singular behaviour at the corner $\bP_j$ we take $\bv_1^j$ to be a piecewise polynomial on a mesh graded geometrically towards $\bP_j$ (more details of these meshes are given below).

A typical term $\bv_2^{j}(\bx) \exp[\ri k_2 r_{j}(\bx)]$ in the second sum in \rf{eqn:ansatztransmission} corresponds to a diffracted wave emanating from the corner $\bP_j$ and propagating with wavenumber $k_2$ in the interior domain $\Omega_2$. By the convexity of $\Omega_2$, this wave impinges on all sides of the polygon. Hence $\bv_2^j(\bx)$ is supported on the whole of $\Gamma$. As for $\bv_1^j$, on each of $\Gamma_j$ and $\Gamma_{j-1}$ we take $\bv_2^j$ to be a piecewise polynomial on a mesh graded geometrically towards $\bP_j$ (again, more details of these meshes are given below).

On $\Gamma\setminus(\Gamma_j\cup\Gamma_{j-1})$ we take $\bv_2^j$ to be a piecewise polynomial of maximum degree $p\in\N_0$ on a mesh of $\Gamma\setminus(\Gamma_j\cup\Gamma_{j-1})$ constructed in the following way. We start with a mesh of $n_s-2$ elements, consisting of the $n_s-2$  sides $\Gamma_l$, $l\not\in\{j,j-1\}$ (with $n_s-1$ mesh points at the corners $\bP_l$, $l\neq j$). We then insert mesh points at any beam boundary points $\bx\in\Gamma\setminus(\Gamma_j\cup\Gamma_{j-1})$ associated with ``strong'' beam boundaries (as defined in \rf{tolBBDef}) in $\bv^{GO}$, provided that the beam boundaries in question emanated from the corner $\bP_j$. (Here we note another difference between the current algorithm and that described in \cite{groth2015hybrid}, where mesh points were introduced at \textit{all} strong beam boundaries, not just those emanating from corner $\mathbf{P}_j$, leading to redundancy in the approximation space due to the same mesh points being added to multiple overlapping meshes. The current algorithm avoids this redundancy, giving a significant reduction in the total number of fewer degrees of freedom, without compromising on accuracy.)

To clarify the mesh refinement procedure, we consider the contribution of the beams illustrated in Figure \ref{fig:rays}. Let us assume that all of the beam boundary discontinuities in Figure \ref{fig:rays} are strong (i.e., satisfy \rf{tolBBDef}). Then for the primary transmitted beam emanating from side $\Gamma_2$ in Figure \ref{fig:rays}(a), we would insert a mesh point at $\bx_5$ in the mesh for $\bv_2^2$ (associated with diffraction from corner $\bP_2$), and at $\bx_4$ in the mesh for $\bv_2^3$ (associated with diffraction from corner $\bP_3$). For the internally reflected beams emanating from side $\Gamma_4$ and $\Gamma_5$ in Figure \ref{fig:rays}(b), we would insert a mesh point at $\bx_3$ in the mesh for $\bv_2^5$ (associated with diffraction from corner $\bP_5$). But we would not insert mesh points at $\bx_2$ and $\bx_6$ in the mesh for any of the amplitudes $\bv_2^j$, $j=1,\ldots,n_s$, as the points $\bx_2$ and $\bx_6$ correspond to beam boundaries that are not directly associated with any corner $\bP_j$. (Rather, they arise by the internal reflection of the two beam boundaries in Figure \ref{fig:rays}(a), and are associated with the higher-order diffracted-reflected fields, which are not captured in our approximation space.) 

To describe in more detail the geometrically graded meshes we use, 
we first consider a mesh on the interval $[0,1]$, refined towards $0$. 
Given $n\geq 1$ (the number of layers in the mesh) we let $G_n$ denote the set of meshpoints $\{x_i\}_{i=0}^n$ defined by
\begin{equation}
      x_0:=0,\quad x_i:=\sigma^{n-i},\quad i=1,2,\ldots,n,  
      \label{mesh_graded}
\end{equation}
where $0<\sigma <1$ is a grading parameter (we discuss the choice of $\sigma$ in \S\ref{sec:implementation}). 
We then define a space of piecewise polynomials on the mesh $G_n$, for a given degree vector $\mathbf{p}\in(\mathbb{N}_0)^n$, by
\[
   \left\{\rho:[0,1]\rightarrow\mathbb{C}:\rho|_{(x_{i-1},x_i)} \text{ is a polynomial of degree less than or equal to } (\mathbf{p})_i,\, i=1,\ldots,n\right\}.
\]
For reasons of efficiency and conditioning it is common to decrease the degree of the approximating polynomials towards the point of refinement. Specifically, given a maximum polynomial degree $p\in\N_0$, we use a ``linear slope'' degree vector $\mathbf{p}$ with
\begin{equation}
   (\mathbf{p})_i:= \left\{ 
  \begin{array}{l l}
    p-\left\lfloor \frac{(n+1-i)}{n}p\right\rfloor, & \quad 1\leq i\leq n-1,\\
    p, & \quad i=n.
  \end{array} \right.
  \label{eqn:poly_deg}
\end{equation} 
We adopt an ``$hp$'' approximation strategy, in which the number of layers $n$ in each geometric mesh grows as the maximum polynomial degree $p$ is increased. Specifically, we take
\begin{align}
\label{CnpDef}
n=\lceil C_{np}(p+1) \rceil,
\end{align}
where $C_{np}>0$ is a user-specified constant. 
The geometric meshes in our approximation space are constructed from this basic building block by straightforward coordinate transformations. Explicitly, for each $j=1,\ldots,n_s$ the amplitudes $\bv_1^j$ and $\bv_2^j$ are approximated on the side $\Gamma_j$ using the mesh 
$\bP_j + G_n*(\bP_{j+1}-\bP_j)$,
and on the side $\Gamma_{j-1}$ using the mesh $\bP_j + G_n*(\bP_{j-1}-\bP_j)$. 

The above construction constrains $\bv^d$ to lie in a particular finite-dimensional approximation space $V\subset L^2(\Gamma)\times L^2(\Gamma)$ 
whose dimension $N$
(the total number of degrees of freedom) is given  by
\begin{equation}
 N = 2\left((p+1)(n_s-2 + n_{bb}) + 4n_s \sum_{i=1}^n ((\bp)_i+1)\right),
 \label{eqn:num_DOF}
 \end{equation}
where $n_{bb}$ is the total number of strong beam boundary points at which we insert extra mesh points.
Note that, by \rf{CnpDef}, $N$ grows approximately in proportion to $p^2$ as $p$ increases.

\subsection{Numerical best approximation error analysis}
\label{sec:NumBestApprox}
The HNA approximation strategy described above has been extensively tested in \cite{GrothPhDThesis} and \cite{groth2015hybrid} using a numerical best approximation error analysis. 
For a number of scattering scenarios, highly accurate reference solutions, denoted here by $\mathbf{v}_{\rm ref}$, were computed using a conventional BEM (see \S\ref{sec:Numerical} for more details). The corresponding GO approximations $\mathbf{v}^{GO}$ were calculated using the BTA, and numerical best approximations to $\mathbf{v}_{\rm ref}-\mathbf{v}^{GO}$ (which we identify with $\mathbf{v}^d$, see (\ref{eqn:vDecomp})) from the HNA approximation space $V$ were then computed using orthogonal projection in $L^2(\Gamma)\times L^2(\Gamma)$ (via the least squares approach detailed in \cite{groth2015hybrid}), and the corresponding errors recorded.

\begin{figure}[t!]
\centering
\subfigure[$\mu=1.5+0.00625\ri$]{
\begin{tikzpicture}[scale=1]
\begin{semilogyaxis}[
	xlabel={Polynomial degree $p$},
	ylabel={Relative $L^2$ error in $u_1$ on $\Gamma$},
	ymax={1e-1},
	ymin={1e-5},
	grid=both,
	font=\footnotesize,
	xtick={-1,0,1,2,3,4,5,6,7},
	xticklabels={GO,0,1,2,3,4,5,6,7},
	legend style={at={(0.05,0.05)},anchor=south west}
]
	
\addplot[color=black,mark=*] coordinates {
	(-1 , 9.67e-2)  (0, 5.69e-2)  (1 , 2.74e-2)  (2 , 1.10e-2)  (3 , 1.06e-2)  (4 , 6.81e-3)  (5, 5.99e-3)  (6, 5.12e-3)  (7 , 4.66e-3)
}; 
	
\addplot[color=red,mark=diamond*] coordinates {
	(-1, 6.72e-2) (0, 4.68e-2)  (1 , 1.79e-2)  (2 , 3.20e-3)  (3 , 2.79e-3)  (4 , 1.57e-3)  (5, 1.43e-3)  (6, 1.26e-3)  (7 , 1.07e-3)
}; 
	
\addplot[color=blue,mark=square*] coordinates {
	(-1, 4.96e-2) (0, 3.82e-2)  (1 , 2.43e-2)  (2 , 3.00e-3)  (3 , 1.05e-3)  (4 , 4.38e-4)  (5, 3.31e-4)  (6, 3.11e-4)  (7 , 3.07e-4)
}; 
\legend{$k_1=40$,$k_1=80$,$k_1=160$}
\end{semilogyaxis}
\end{tikzpicture}
}
\subfigure[$\mu=1.5+0.0125\ri$]{
\begin{tikzpicture}[scale=1]
\begin{semilogyaxis}[
	xlabel={Polynomial degree $p$},
		ymax={1e-1},
	ymin={1e-5},
	grid=both,
	font=\footnotesize,
	xtick={-1,0,1,2,3,4,5,6,7},
	xticklabels={GO,0,1,2,3,4,5,6,7},
]
	
\addplot[color=black,mark=*] coordinates {
	(-1 , 9.14e-2) (0, 5.28e-2)  (1 , 1.54e-2)  (2 , 3.39e-3)  (3 , 3.00e-3)  (4 , 1.48e-3)  (5, 1.30e-3)  (6, 1.10e-3)  (7 , 9.99e-4)
}; 
	
\addplot[color=red,mark=diamond*] coordinates {
	(-1 , 6.77e-2) (0, 4.59e-2)  (1 , 1.76e-2)  (2 , 1.91e-3)  (3 , 1.23e-3)  (4 , 5.10e-4)  (5, 3.55e-4)  (6, 3.24e-4)  (7 , 3.14e-4)
}; 
	
\addplot[color=blue,mark=square*] coordinates {
	(-1 , 4.93e-2)  (0, 3.58e-2)  (1 , 2.42e-2)  (2 , 2.86e-3)  (3 , 8.27e-4)  (4 , 2.94e-4)  (5, 9.29e-5)  (6, 3.80e-5)  (7 , 2.12e-5)
}; 
 
\end{semilogyaxis}
\end{tikzpicture}
}
\caption{Numerical best approximation errors for the configuration in Figure \ref{fig:tri_decomp} 
at two different refractive indices (left panel - lower absorption, right panel - higher absorption), as a function of polynomial degree $p$. Other approximation space parameters ($C_{np}$ etc.) are as specified in \S\ref{sec:Numerical}. For reference the GO error is also shown.}
\label{fig:convergence_p}
\end{figure}
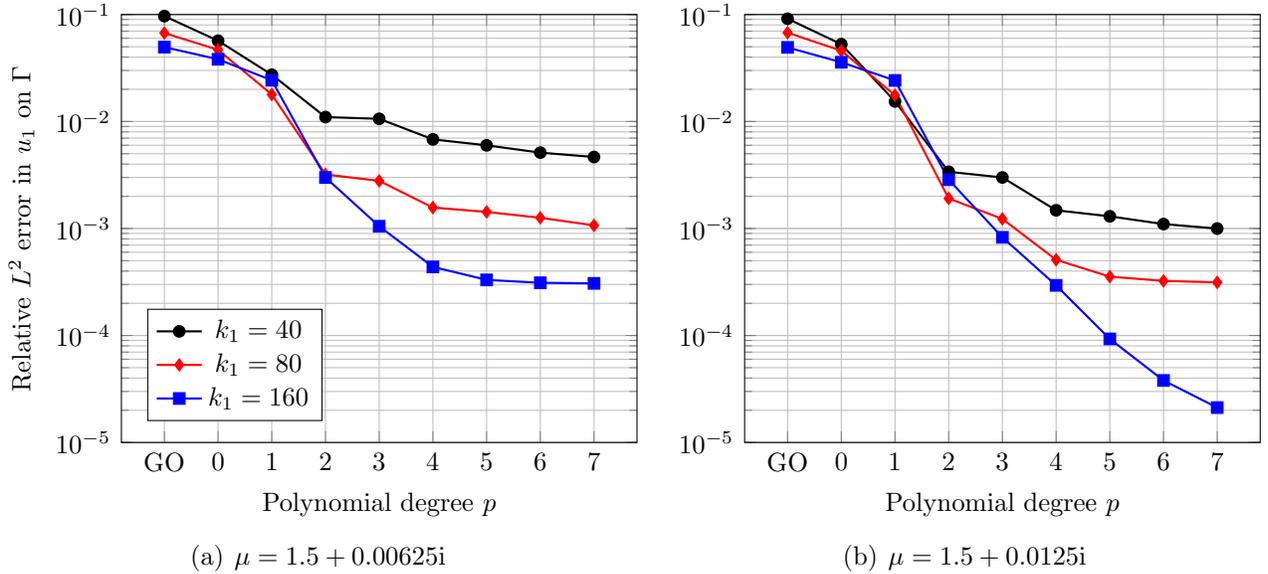

A typical sample of the results obtained is presented in Figure \ref{fig:convergence_p}, which plots the numerical best approximation error against increasing polynomial degree. For reference we also show the error for GO alone, as the left-most data point. All of the curves in Figure \ref{fig:convergence_p} exhibit the same qualitative behaviour. At first the error decays exponentially, in an essentially frequency-independent way; in this initial phase, increasing $p$ leads to an improvement in the accuracy of our HNA approximation $\bv^d$ to the diffracted wave fields. Then, when the accuracy of $\bv^d$ reaches the level of the asymptotic error associated with our neglect of the higher-order asymptotic phenomena (the lateral waves and the internal re-reflections of the diffracted and lateral waves), the error levels out. Increasing $p$ no longer leads to a significant improvement in accuracy, because our HNA approximation space for $\bv^d$ does not contain the correct oscillatory basis functions needed to capture the neglected waves. (Exponential convergence will in theory resume once the number of degrees of freedom is sufficiently large to resolve the missing oscillations, but this is not the regime in which HNA methods are designed to operate.) 
We note that the magnitude of the neglected waves, and hence the level at which the error stagnates, decreases as $k_1$ increases, and also as the absorption $\im{\mu}$ increases, in accordance with the asymptotic theory (see the discussion before \S\ref{sec:comp_v_go}).

\section{Galerkin BEM and numerical results}
\label{sec:Numerical}

In \S\ref{sec:HNA} we proposed an ansatz $\bv\approx \bv^{GO}+\bv^d$ (see (\ref{eqn:vDecomp})) for the solution $\bv=(u_1,\pdonetext{u_1}{\bn})^T\in L^2(\Gamma)\times L^2(\Gamma) $ of the integral equation (\ref{matrix}), 
describing a beam tracing algorithm for computing the GO term $\bv^{GO}$ and an HNA approximation space $V\subset L^2(\Gamma)\times L^2(\Gamma)$ for the diffraction term $\bv^d$. 
In this section we describe our BEM implementation and present a range of numerical results demonstrating its performance. Our BEM selects a particular element $\bv^d\in V$ by applying the Galerkin method to the integral equation~(\ref{matrix}), rewritten using the decomposition~(\ref{eqn:vDecomp}), with $\bv^d$ as the unknown.  That is, we compute $\bv^d\in V$ satisfying
\begin{equation}
  \langle A\bv^d, \bw \rangle = \langle \mathbf{f} - A\bv^{GO}, \bw \rangle, \quad \mbox{for all } \bw \in V,
  \label{eqn:gal}
\end{equation}
where $\langle \cdot , \cdot\rangle$ is the usual inner product on $L^2(\Gamma)\times L^2(\Gamma)$.

\subsection{Implementation details}
\label{sec:implementation}
Choosing a basis for $V$ reduces \rf{eqn:gal} to the solution of a linear system for the basis coefficients of $\bv^d$. In our implementation we represent each of the two components of $\bv^d$ using scalar basis functions of the form $L_m \re^{\ri k_l r_j}$ (see \rf{eqn:ansatztransmission}), where $l\in\{1,2\}$, $r_j$ is the distance from vertex $\bP_j$, for some $j\in\{1,\ldots,n_s\}$, and $L_m$ is a Legendre polynomial of degree $m\in\{0,\ldots,p\}$. (Note that each basis function is supported on some subinterval of one of the sides of the polygon.) 
To improve conditioning of the resulting linear system we found it beneficial to scale each basis function so that its $L^2(\Gamma)$ norm is equal to one. 

By design, the linear systems produced by the HNA BEM are small compared to those arising from conventional methods 
and hence can be solved using direct methods, without the need for iterative solvers. Assembly of the linear system requires the evaluation of integrals which are potentially both oscillatory and singular. We evaluate these integrals to high precision (better than $10^{-8}$ relative accuracy) with composite Gauss quadrature using geometric grading and the classical Duffy transformation - for details see \cite{GrothPhDThesis}. For a faster implementation one could employ oscillatory quadrature techniques, as discussed in \cite{acta,DoGrKi:13,kim2012asymptotic,hewett2014frequency,hewett2013high}. 
But our focus in the current paper is on demonstrating the efficiency of the HNA approximation space in terms of the number of degrees of freedom, 
rather than on fast implementations, 
and so we postpone further discussion of this important issue to future work.

In all of our experiments we use the parameter choices:
\begin{align}
\label{tol_b}{\rm tol}_{b}=0.005 \qquad \textrm{(see \rf{tolBDef})},\\
\label{tol_GO}{\rm tol}_{GO}=0.01 \qquad \textrm{(see \rf{tolGODef})},\\
\label{tol_BB} {\rm tol}_{bb}=0.01 \qquad \textrm{(see \rf{tolBBDef})},\\
\label{C_np}C_{np}=1.5 \qquad \textrm{(see \rf{CnpDef})},
\end{align}
and the fixed maximum polynomial degree
\begin{align}
\label{p}
p=3 \qquad \textrm{(see \rf{eqn:poly_deg})}.
\end{align}
In accordance with the discussion in \S\ref{sec:NumBestApprox}, choosing a larger value of $p$ than $3$ in general leads to improved accuracy, up to the point at which the underlying asymptotic error is reached. However, through extensive numerical testing (detailed in \cite{GrothPhDThesis}) the choices \rf{tol_b}--\rf{p} were found to give a reasonable compromise between computational effort and solution accuracy for the range of problems considered here.

For the construction of our geometric meshes (recall~(\ref{mesh_graded})) we use grading parameters $\sigma_1=0.17$ for the $\bv_1^j$ meshes and $\sigma_2=0.15$ for the $\bv_2^j$ meshes. 
Using the same grading parameter for both the $\bv_1^j$ and $\bv_2^j$ meshes was found to cause ill-conditioning when $k_1$ is small, or when the meshes are particularly heavily refined, because when the diameter of the smallest elements in the meshes is comparable to the wavelength(s), the polynomial approximants are capable of resolving the oscillatory difference between the two phase factors $\exp[\ri k_1r_j]$ and $\exp[\ri k_2r_j]$, leading to redundancy in the approximation space. 
Of course our method is designed with the high frequency regime in mind, and this redundancy is a low-frequency issue. But using slightly different grading parameters for the two meshes was found to be a simple and effective way to maintain stability even for relatively low frequencies.

\subsection{Numerical experiments}
\label{subsec:num}
In \S\ref{sec:vary_freq}--\S\ref{subsec:geom} below, we present numerical results demonstrating the accuracy and efficiency of our HNA BEM for a range of scattering scenarios. As we shall see, our HNA BEM achieves an order of magnitude improvement in accuracy over the purely asymptotic PGOH approach, using only a modest (and frequency-independent) number of degrees of freedom. 

To quantify the accuracy of the PGOH and HNA BEM solutions we compare them to reference solutions obtained by solving (\ref{matrix}) using a conventional $hp$-BEM with a large number of degrees of freedom. Numerical experiments suggest that this reference solution achieves at least $0.01\%$ accuracy for all the problems considered here; for more details see \cite{GrothPhDThesis}. 
We shall present results both for the boundary solution $\bv=(u_1,\pdonetext{u_1}{\bn})^T$ and the far-field pattern $F$ (see \eqref{eqn:2D_FF}).

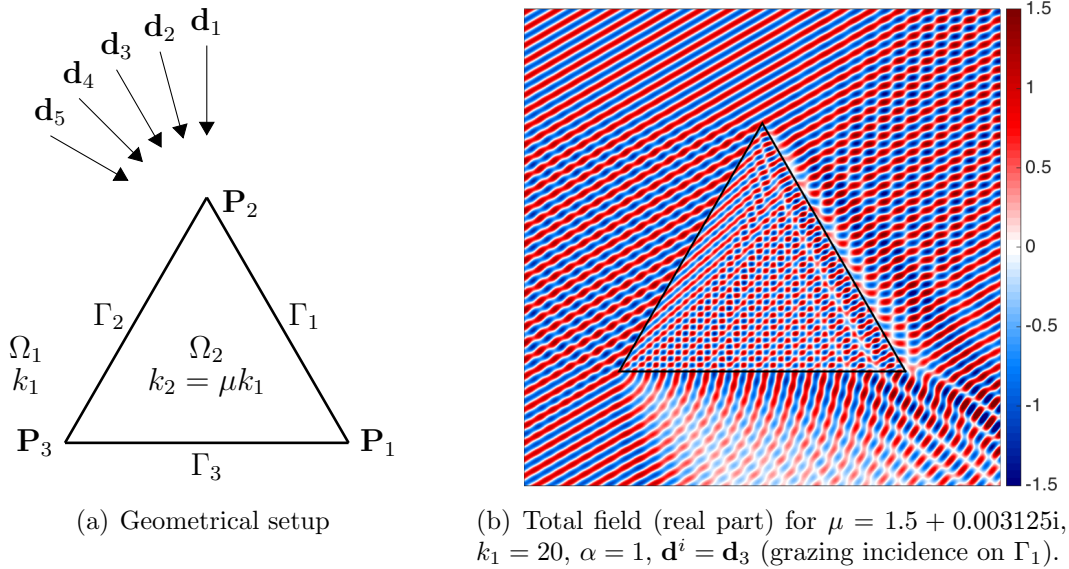
\begin{figure}[t]
\centering
\subfigure[Geometrical setup]
{
\begin{tikzpicture}[scale=0.6]
    \draw[line width=1pt](-pi,-1.81)--(pi,-1.81)node[anchor=west]{$\mathbf{P}_1$};
    \draw[line width=1pt](pi,-1.81)--(0,3.63);
    \draw(0.75,3.5) node{$\mathbf{P}_2$};
    \draw[line width=1pt](0,3.63)--(-pi,-1.81)node[anchor=east]{$\mathbf{P}_3$};             
    \draw[arrows={-triangle 60}];
        \path(0,-1.81)node[anchor=north]{$\Gamma_3$};
    \path(-pi/2,1)node[anchor=east]{$\Gamma_2$};
    \path(pi/2,1)node[anchor=west]{$\Gamma_1$};
        \draw(0,0.8)node[anchor=north]{$\Omega_2$};
                \draw(0,0.1)node[anchor=north]{$k_2=\mu k_1$};
           \draw(-4,0.8)node[anchor=north]{$\Omega_1$};
           \draw(-4,0.1)node[anchor=north]{$k_1$};
\begin{scope}[xshift=0cm,yshift=3cm]
\draw[arrows={-triangle 60}] (90:4)node[anchor=south]{$\bd_1$}--(90:2);
\draw[arrows={-triangle 60}] (105:4)node[anchor=south]{$\bd_2$}--(105:2);
\draw[arrows={-triangle 60}] (120:4)node[anchor=south]{$\bd_3$}--(120:2);
\draw[arrows={-triangle 60}] (135:4)node[anchor=south]{$\bd_4$}--(135:2);
\draw[arrows={-triangle 60}] (150:4)node[anchor=south]{$\bd_5$}--(150:2);
\end{scope}    
\end{tikzpicture}
}
\hs{5}
\subfigure[Total field (real part) for $\mu=1.5+0.003125\ri$, $k_1=20$, $\alpha=1$, $\bd^i=\bd_3$ (grazing incidence on $\Gamma_1$).]{\hs{5}
\includegraphics[width=0.4\textwidth]{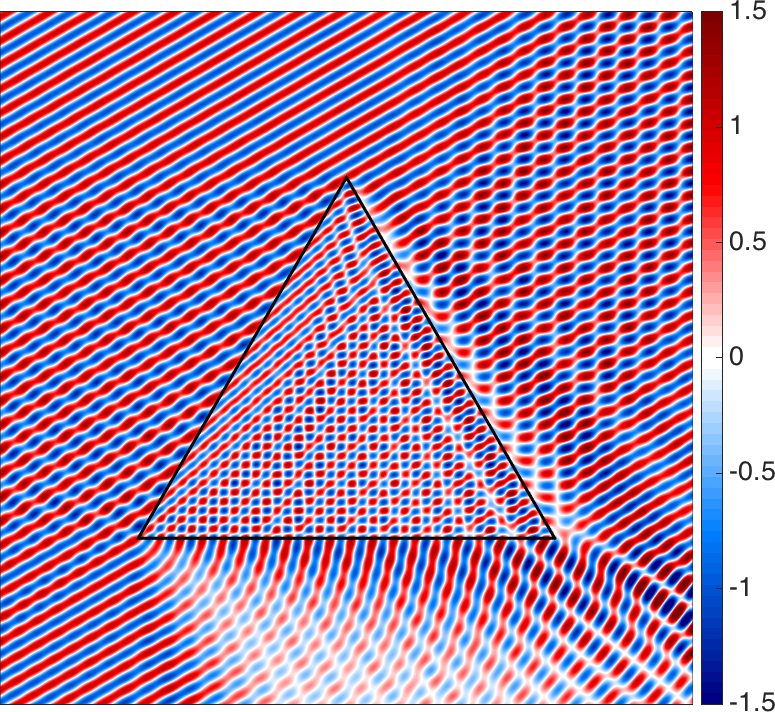}
}
\caption{Scattering by an equilateral triangle. }
\label{fig:triangle_setup}
\end{figure}

Our algorithm can be applied to any convex polygonal scatterer, and to illustrate this we consider three different examples: an equilateral triangle, a square, and a regular hexagon. In each case, we take the sides of the polygon to be of length $2\pi$ so that the number of exterior (interior) wavelengths $\lambda_1$ ($\lambda_2$) around $\Gamma$ is equal to $n_s k_1$ ($n_s\real{k_2}=n_sk_1\real{\mu}$). We shall report results for $k_1$ ranging from $5$ to $160$, the highest value corresponding to $480$, $640$ and $960$ exterior wavelengths around $\Gamma$ for the triangle, square, and hexagon, respectively. (The reason we do not present results for even higher wavenumbers is that computing the reference solution becomes prohibitively expensive.) 
For comparison with conventional methods (such as the standard $hp$-BEM used to compute our reference solution) 
we introduce the notation
\[ \DOFone = N/(2n_sk_1) \qquad \text{ and } \qquad  \DOFtwo = N/(2n_s\real{k_2})=(\DOFone)/\real{\mu}\]
to denote the average number of degrees of freedom (DOF) per wavelength around $\Gamma$ used in the approximation of each of the components of the boundary solution $\bv$, relative to the exterior and interior wavelengths respectively. 
A rule of thumb for scattering problems (see, e.g., \cite{sixDOF,marburg2008discretization}) is that to achieve acceptable ``engineering accuracy'' using a conventional method, one generally requires 6 to 10 degrees of freedom per wavelength. Precisely what is meant by ``engineering accuracy'' depends on the application at hand, but in \cite{sixDOF,marburg2008discretization} the target accuracy appears to be roughly $10\%$.  As we shall see, at high frequencies our HNA BEM can achieve better than $1\%$ accuracy with far less than one degree of freedom per wavelength.

For the equilateral triangle we shall consider scattering for the five different incident directions $\bd_1$, $\bd_2$, $\bd_3$, $\bd_4$ and $\bd_5$ shown in Figure \ref{fig:triangle_setup}(a), 
which correspond to taking 
\begin{align}
\label{eqn:diDef}
\bd^i = (\cos\theta_i,-\sin\theta_i), \quad \text{with } \theta_i = \frac{\pi}{2}, \frac{5\pi}{12}, \frac{\pi}{3}, \frac{\pi}{4}, \frac{\pi}{6},
\end{align}
respectively. 
We also study the performance of our method under changes in the refractive index $\mu$, and the constant $\alpha$ in (\ref{bcs}).  
Our work is motivated in part by applications in atmospheric physics, specifically the scattering of electromagnetic radiation by ice crystals in cirrus clouds \cite{Baran:2009}. In this context the relevant choices of $\alpha$ are $\alpha=1$ and $\alpha=(k_1/k_2)^2$, as already discussed after \rf{SRC}. The refractive index $\mu$ of ice has been determined through numerous experiments (see, e.g., \cite{warren2008optical}) and is known to be complex-valued and highly frequency-dependent, with 
$\real{\mu}$ 
(the \emph{contrast}) 
ranging between $1$ and $2$, and 
$\im{\mu}$ 
(the \emph{absorption}) 
between $0$ and $1$. 
For most of our examples we take $\real{\mu}=1.5$. 
But we also present results for $\real{\mu}=0.66$ (the approximate reciprocal of $1.5$), to demonstrate that our method is also effective for the case $0<\real{\mu}<1$. 
We consider a range of $\im{\mu}$ from $0$ up to $0.0125$; results for larger absorptions (when our method achieves accuracies even better than those reported here) can be found in \cite{GrothPhDThesis}. 

\subsection{High frequency performance}
\label{sec:vary_freq}
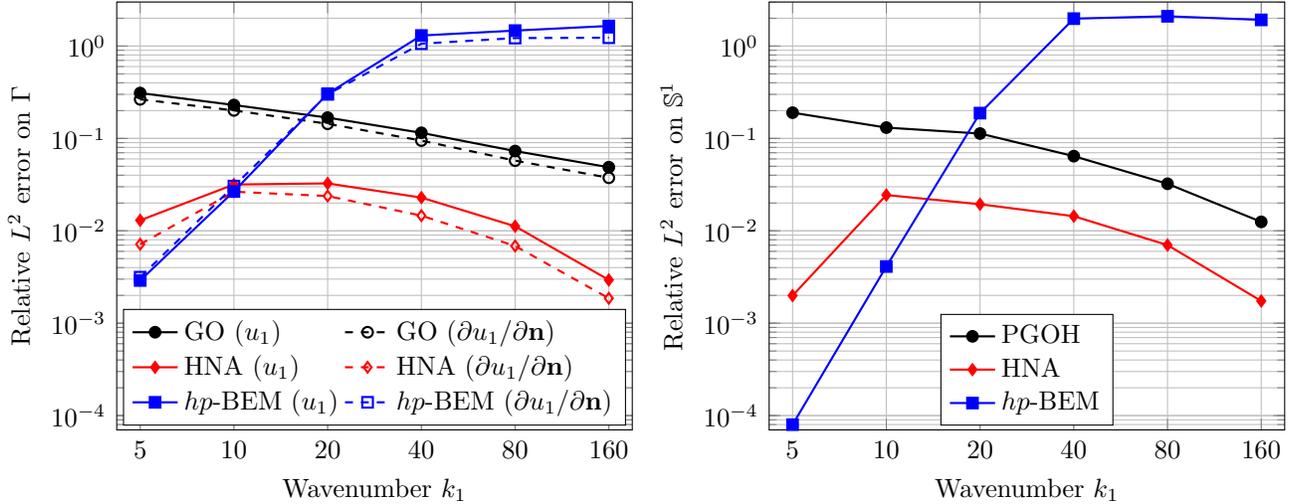
\begin{figure}[t!]
\centering
\subfigure[Boundary solution $\bv=(u_1,\pdonetext{u_1}{\bn})^T$]{
\begin{tikzpicture}[scale=1]
\begin{loglogaxis}[
	xlabel={Wavenumber $k_1$},
	ylabel={Relative $L^2$ error on $\Gamma$},
	ymax={3e0},
	ymin={7e-5},
	grid=both,
	xtick={5,10,20,40,80,160},
	xticklabels={5,10,20,40,80,160},
	legend style={at={(0.5,0.01)},anchor=south},
	legend columns=2,
	legend cell align={left},
	font=\footnotesize,
	enlarge x limits=0.05
]

\addplot[color=black,mark=*] coordinates {
	(5, 3.10e-1) (10, 2.30e-1) (20, 1.68e-1) (40, 1.15e-1) (80, 7.30e-2) (160, 4.88e-2)
}; 
\addlegendentry{GO ($u_1$)}
\addplot[color=black,dashed,mark=o,mark options={solid}] coordinates {
	(5, 2.65e-1) (10, 2.01e-1) (20, 1.44e-1) (40, 9.48e-2) (80, 5.75e-2)  (160, 3.75e-2)
}; 
\addlegendentry{GO ($\pdonetext{u_1}{\bn}$)}
	
\addplot[color=red,mark=diamond*] coordinates {
	(5,1.30e-2) (10, 3.17e-2) (20, 3.26e-2) (40, 2.29e-2) (80,1.12e-2) (160,2.94e-3)
}; 
\addlegendentry{HNA ($u_1$)}

\addplot[color=red,dashed,mark=diamond,mark options={solid}] coordinates {
	(5, 7.15e-3) (10, 2.65e-2) (20, 2.38e-2) (40, 1.46e-2) (80, 6.86e-3) (160, 1.86e-3)
}; 
\addlegendentry{HNA ($\pdonetext{u_1}{\bn}$)}	

\addplot[color=blue,mark=square*] coordinates {
	(5,2.90e-3) (10, 2.67e-2) (20, 3.05e-1) (40, 1.30e0) (80,1.47e0) (160, 1.65e0)
}; 
\addlegendentry{$hp$-BEM ($u_1$)}

\addplot[color=blue,dashed,mark=square,mark options={solid}] coordinates {
	(5, 3.15e-3) (10, 3.05e-2) (20, 3.00e-1) (40, 1.06e0) (80,1.22e0) (160, 1.23e0)
}; 	
\addlegendentry{$hp$-BEM ($\pdonetext{u_1}{\bn}$)}	
\end{loglogaxis}
\end{tikzpicture}
}
\hs{-4}
\subfigure[Far-field pattern $F$]{
\begin{tikzpicture}[scale=1]
\begin{loglogaxis}[
	xlabel={Wavenumber $k_1$},
	ylabel={Relative $L^2$ error on $\mathbb{S}^1$},
	ymax={3e0},
	ymin={7e-5},
	grid=both,
	xtick={5,10,20,40,80,160},
	xticklabels={5,10,20,40,80,160},
	legend style={at={(0.5,0.01)},anchor=south},
	legend cell align={left},
	font=\footnotesize,
	enlarge x limits=0.05
]

\addplot[color=black,mark=*] coordinates {
	(5, 1.90e-1) (10, 1.31e-1) (20, 1.13e-1) (40, 6.43e-2) (80, 3.23e-2)  (160, 1.25e-2)
}; 
\addlegendentry{PGOH}	
	
\addplot[color=red,mark=diamond*] coordinates {
	(5,1.99e-3) (10, 2.44e-2) (20, 1.94e-2) (40, 1.44e-2) (80, 7.01e-3) (160, 1.74e-3)
}; 
\addlegendentry{HNA}	
	
\addplot[color=blue,mark=square*] coordinates {
	(5,7.98e-5) (10, 4.11e-3) (20, 1.88e-1) (40, 1.98e0) (80, 2.10e0) (160, 1.92e0)
}; 
\addlegendentry{$hp$-BEM}		
\end{loglogaxis}
\end{tikzpicture}
}
\caption{Scattering by an equilateral triangle with $\bd^i=\bd_1$, $\alpha=1$ and $\mu=1.5+0.003125\ri$.  Errors in boundary solution and far-field pattern for GO (PGOH in the far field), our HNA BEM (with fixed $N=416$) and a conventional $hp$-BEM (with fixed $N=456$). 
}
\label{fig:fixed_DOF_456}
\end{figure}

In this section we investigate how the accuracy of our method depends on the wavenumber $k_1$. 
We consider the configuration from Figure~\ref{fig:tri_decomp}, namely scattering by an equilateral triangle with $\bd^i=\bd_1$, $\alpha=1$ and $\mu=1.5+0.003125\ri$, for wavenumbers $k_1=5,10,20,40,80,160$. The parameter choices \rf{tol_b}--\rf{p} result in a total number of degrees of freedom $N=416$ in this case, independent of $k_1$. 

Figure \ref{fig:fixed_DOF_456} shows the resulting error in our HNA BEM approximation to the boundary solution $\bv=(u_1,\pdonetext{u_1}{\bn})^T$ and the far-field pattern $F$. 
For a comparison with a purely asymptotic method we also show the corresponding errors for the GO approximation, and the resulting PGOH approximation to the far-field pattern (obtained by substituting the GO approximation for $\bv$ into \eqref{eqn:2D_FF}). 
For a comparison with a purely numerical method we also show the corresponding errors for a ``conventional'' BEM, namely the same $hp$-BEM used to generate our reference solution, but with far fewer degrees of freedom; specifically, 
for the $hp$-BEM we use approximately the same ($k_1$-independent) number of degrees of freedom (in this case $N=456$) 
as we use in the HNA BEM ($N=416$). 
The corresponding values of $\DOFtwo$ for the two methods, along with the 2-norm condition numbers COND of the associated linear systems, are presented in Table \ref{tab:DOFperlamAndCond}. 

For all values of $k_1$ considered here, the HNA method is accurate to within at least $5\%$ both for the computation of the boundary solution and the far-field pattern, and provides an order of magnitude improvement over GO/PGOH. Moreover, as is true for GO/PGOH, the accuracy of the HNA method \emph{increases} as $k_1$ increases, despite the fact that $N$ remains fixed. For the largest wavenumber $k_1=160$, the HNA method achieves an accuracy of approximately $0.1\%$ using fewer than $0.3$ degrees of freedom per wavelength. In contrast, the error for the conventional $hp$-BEM (again with fixed $N$), despite starting below that of the HNA method for $k_1=5$, rapidly increases with increasing $k_1$, with the accuracy cross-over between $hp$-BEM and HNA BEM occurring between $k_1=10$ and $k_1=20$. The $hp$-BEM loses accuracy completely (relative error $>100\%$) for $k_1\geq 40$. (Of course, accuracy could be regained for the $hp$-BEM for $k_1\geq 40$ by increasing the number of DOF, but the point of Figure \ref{fig:fixed_DOF_456} is to compare the performance of the conventional $hp$-BEM and our HNA BEM when they are allocated a \emph{fixed} (and roughly equal) number of DOFs.) 

Table~\ref{tab:DOFperlamAndCond} reveals that, for the examples considered, the condition number for the HNA BEM broadly decreases as frequency increases, in contrast to the conventional $hp$-BEM.  This is consistent with the fact that as $k_1$ increases our HNA approximation space becomes increasingly well-suited to approximating the boundary solution.  The conditioning presents no issues for the direct linear solver used in the HNA BEM, as the error plots in Figure~\ref{fig:fixed_DOF_456} demonstrate. 

\begin{table}[t]
\centering
\begin{tabular}{|l|l|l|l|l|l|l|}
\hline
$k_1$ & 5 & 10 & 20 & 40 & 80 & 160 \\
\hline
$hp$-BEM $\DOFtwo$ & 10.13 & 5.07 & 2.53 & 1.27 & 0.63 & 0.32\\
\hline
HNA $\DOFtwo$ & 9.24 & 4.62 & 2.31 & 1.16 & 0.58 & 0.29\\
\hline
$hp$-BEM COND & $6.95\times 10^1$ & $2.54\times 10^2$  & $5.51\times 10^2$ & $1.53\times 10^3$ & $4.94\times 10^3$ & $3.20\times 10^4$ \\
\hline
HNA BEM COND & $1.42\times 10^8$  & $1.55\times 10^6$ &  $1.93\times 10^6$ & $4.77\times 10^5$ & $1.46\times 10^5$ & $1.08\times 10^5$ \\
\hline	
\end{tabular}
\caption{Degrees of freedom per wavelength $\lambda_2$ for the results in Figure \ref{fig:fixed_DOF_456}. The total number of DOF is 456 for the $hp$-BEM and 416 for the HNA BEM, independent of $k_1$.
\label{tab:DOFperlamAndCond}
}
\end{table}

\subsection{Varying the incident angle}
\label{sec:vary_inc_angle}
We now investigate how the accuracy of the method depends on the incident wave direction. As in the last section we fix $\alpha=1$ and $\mu = 1.5+0.003125\ri$, but now consider the five evenly spaced incident angles $\bd_1,\ldots,\bd_5$ defined by \rf{eqn:diDef} (see Figure \ref{fig:triangle_setup}(a)). 
We already presented field plots for the case $\bd^i=\bd_1$ in Figure~\ref{fig:tri_decomp}. Figure \ref{fig:triangle_setup}(b) shows an analogous field plot for the case $\bd^i=\bd_3$, which corresponds to grazing incidence along $\Gamma_1$. 
The number of degrees of freedom $N$ used in the approximation space for each of these examples ranges from 400 and 424; exact values are detailed in Table~\ref{tab:DOF_v_angle}. The variation between examples is due to the variation of the strength of beam boundary discontinuities, and whether or not extra mesh points are introduced to capture them (recall \rf{eqn:num_DOF}, \rf{tolBBDef} and the discussion at the end of \S\ref{sec:comp_v_go}).

\begin{table}[hb!]
\centering
\begin{tabular}{| c || c | c | c | c | c | }
\hline 
$\theta_i$  \textbackslash $k_1$ & 10 & 20 & 40 & 80 & 160 \\
\hline\hline 
$\frac{\pi}{2}$ ($\bd_1$) & 416 & 416 & 416 & 416 & 416 \\ \hline
$\frac{5\pi}{12}$ ($\bd_2$) & 424 & 424 & 424 & 416 & 408 \\ \hline
$\frac{\pi}{3}$ ($\bd_3$) & 408 & 408 & 408 & 400 & 400 \\ \hline
$\frac{\pi}{4}$ ($\bd_4$) & 416 & 408 & 408 & 408 & 400 \\ \hline
$\frac{\pi}{6}$ ($\bd_5$) & 408 & 408 & 408 & 408 & 400 \\ \hline
	
\end{tabular}
\caption{Number of DOF in the HNA approximation space for the examples considered in Figure~\ref{fig:vary_angle}.}
\label{tab:DOF_v_angle}
\end{table}

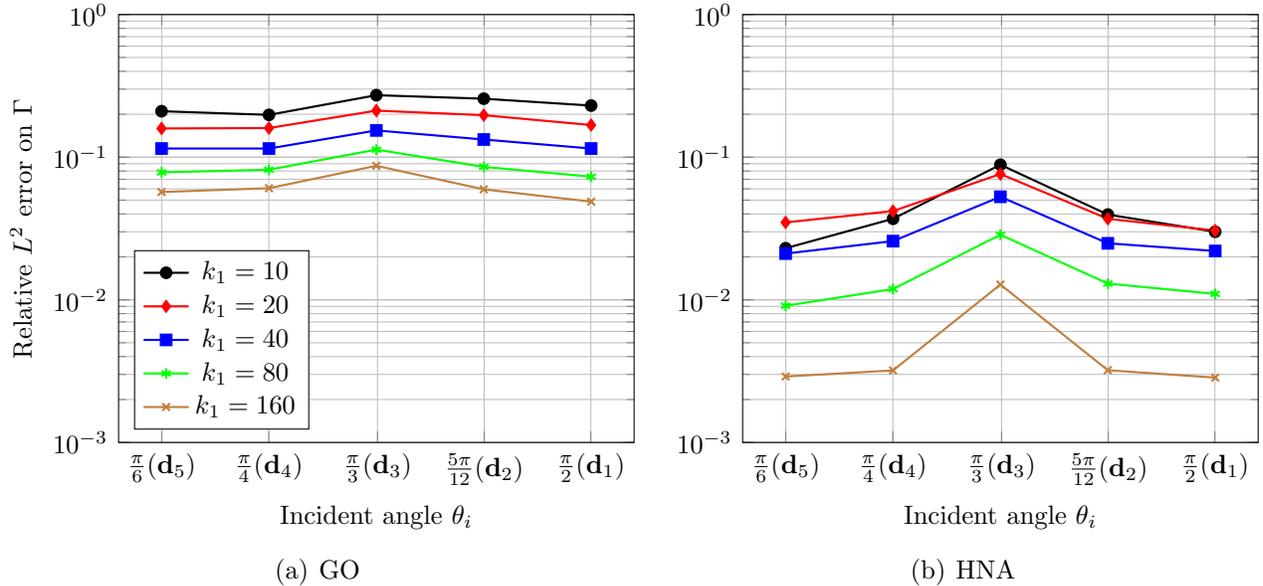
\begin{figure}[t!]
\centering
\subfigure[GO]{
\begin{tikzpicture}[scale=1]
\begin{semilogyaxis}[
	xlabel={Incident angle $\theta_i$},
	ylabel={Relative $L^2$ error on $\Gamma$},
	ymax={1*10^0},
	ymin={1*10^(-3)},
	grid=both,
	xtick={0.524,0.785,1.05,1.31,1.57},
	xticklabels={$\frac{\pi}{6}(\bd_5)$, $\frac{\pi}{4}(\bd_4)$, $\frac{\pi}{3}(\bd_3)$, $\frac{5\pi}{12}(\bd_2)$, $\frac{\pi}{2}(\bd_1)$},
	legend pos={south west},
	font=\footnotesize,
]
	
\addplot[color=black,mark=*] coordinates {
	(pi/6 , 2.10e-1)  (pi/4 , 1.98e-1)  (pi/3 , 2.72e-1)  (5*pi/12 , 2.57e-1)  (pi/2 , 2.30e-1)
};
\addlegendentry{$k_1=10$}
	
\addplot[color=red,mark=diamond*] coordinates {
	(pi/6 , 1.59e-1)  (pi/4 , 1.60e-1)  (pi/3 , 2.12e-1)  (5*pi/12 , 1.97e-1)  (pi/2 , 1.68e-1)
};    
\addlegendentry{$k_1=20$}
	
\addplot[color=blue,mark=square*] coordinates {
	(pi/6 , 1.15e-1)  (pi/4 , 1.15e-1)  (pi/3 , 1.54e-1)  (5*pi/12 , 1.33e-1)  (pi/2 , 1.15e-1)
};      
\addlegendentry{$k_1=40$}
	
\addplot[color=green,mark=asterisk] coordinates {
	(pi/6 , 7.83e-2)  (pi/4 , 8.166e-2)  (pi/3 , 1.13e-1)  (5*pi/12 , 8.55e-2)  (pi/2 , 7.30e-2)
};   
\addlegendentry{$k_1=80$}
	
\addplot[color=brown,mark=x] coordinates {
	(pi/6 , 5.70e-2)  (pi/4 , 6.06e-2)  (pi/3 , 8.70e-2)  (5*pi/12 , 5.95e-2)  (pi/2 , 4.88e-2)
}; 
\addlegendentry{$k_1=160$}
	
\end{semilogyaxis}
\end{tikzpicture}
}
\subfigure[HNA]{
\begin{tikzpicture}[scale=1]
\begin{semilogyaxis}[
	xlabel={Incident angle $\theta_i$},
	ymax={1*10^0},
	ymin={1e-3},
	grid=both,
	legend columns={2},
	font=\footnotesize,
	xtick={0.524,0.785,1.05,1.31,1.57},
	xticklabels={$\frac{\pi}{6}(\bd_5)$, $\frac{\pi}{4}(\bd_4)$, $\frac{\pi}{3}(\bd_3)$, $\frac{5\pi}{12}(\bd_2)$, $\frac{\pi}{2}(\bd_1)$},
]
	
\addplot[color=black,mark=*] coordinates {
	(pi/6 , 2.3e-2)  (pi/4 , 3.7e-2)  (pi/3 , 8.85e-2)  (5*pi/12 , 3.96e-2)  (pi/2 , 3.00e-2)
};
	
\addplot[color=red,mark=diamond*] coordinates {
	(pi/6 , 3.49e-2)  (pi/4 , 4.19e-2)  (pi/3 , 7.62e-2)  (5*pi/12 , 3.70e-2)  (pi/2 , 3.06e-2)
};    
	
\addplot[color=blue,mark=square*] coordinates {
	(pi/6 , 2.11e-2)  (pi/4 , 2.58e-2)  (pi/3 , 5.28e-2)  (5*pi/12 , 2.49e-2)  (pi/2 , 2.20e-2)
};      
	
\addplot[color=green,mark=asterisk] coordinates {
	(pi/6 , 9.09e-3)  (pi/4 , 1.19e-2)  (pi/3 , 2.86e-2)  (5*pi/12 , 1.30e-2)  (pi/2 , 1.10e-2)
};   
	
\addplot[color=brown,mark=x] coordinates {
	(pi/6 , 2.90e-3)  (pi/4 , 3.20e-3)  (pi/3 , 1.28e-2)  (5*pi/12 , 3.21e-3)  (pi/2 , 2.85e-3)
};

\end{semilogyaxis}
\end{tikzpicture}
}

\caption{Accuracy of GO and HNA approximations to $u_1$ on $\Gamma$ plotted against incident angle $\theta_i$ for the incident directions $\bd_1,\ldots,\bd_5$ of Figure \ref{fig:triangle_setup}(a), for $\alpha=1$ and $\mu = 1.5+0.003125\ri$.
}
\label{fig:vary_angle}
\end{figure}

In Figure~\ref{fig:vary_angle}  
we report $L^2$ errors in the approximations to $u_1$ on $\Gamma$ for all five incident directions and a range of wavenumbers $k_1$. (Errors in $\pdonetext{u_1}{\bn}$ and $F$ follow very similar trends but are not reported here.) 
In all cases considered, the HNA BEM provides a significant improvement over GO. 
It is interesting to note, however, that the errors in both the GO and HNA approximations are not uniform across incident angle; rather one observes a peak in the error curves at the grazing incidence case $\bd^i=\bd_3$ (i.e.\ $\theta_i=\pi/3$). 
We hypothesize that the larger error in the GO at grazing incidence is due to the fact that the high frequency asymptotic behaviour in this case involves a particularly prominent diffracted wave along $\Gamma_1$, which is not captured by the GO approximation. (By contrast, the HNA approximation does capture this diffracted wave.)

We hypothesize that the larger error in the HNA approximation at grazing incidence is due to a particularly prominent lateral wave contribution on $\Gamma_3$, which isn't captured by our HNA approximation space. As was suggested in the caption to Figure \ref{fig:LateralAndDiffRef}, it is plausible that grazing incidence might lead to particularly large lateral wave contributions, because of the following argument: while the lateral wave usually compensates for the phase mismatch between the exterior and interior diffracted waves, at grazing incidence one of the lateral waves (in this case the one associated with $\Gamma_1$, since for $\bd^i=\bd_3$ the direction of the incident wave is parallel to that side) has to match up to the incident wave; since the incident wave has lower asymptotic order than the diffracted wave, one might therefore expect the lateral wave to be more prominent in this case.

Lateral waves only propagate inside $\Omega_2$ when $\real{\mu}>1$, which is the case for the examples considered up to now. When $0<\real{\mu}<1$ the lateral waves propagate away from the scatterer in $\Omega_1$, and therefore should not affect the boundary solution. Hence we can test our hypothesis about the role that lateral waves play in the deterioration of accuracy at grazing incidence for the case $\real{\mu}>1$ by studying the case $0<\real{\mu}<1$ and observing whether the same deterioration in accuracy occurs. This we do in the next section.

\subsection{Varying the contrast ($\real{\mu}$)}
\label{sec:vary_contrast}
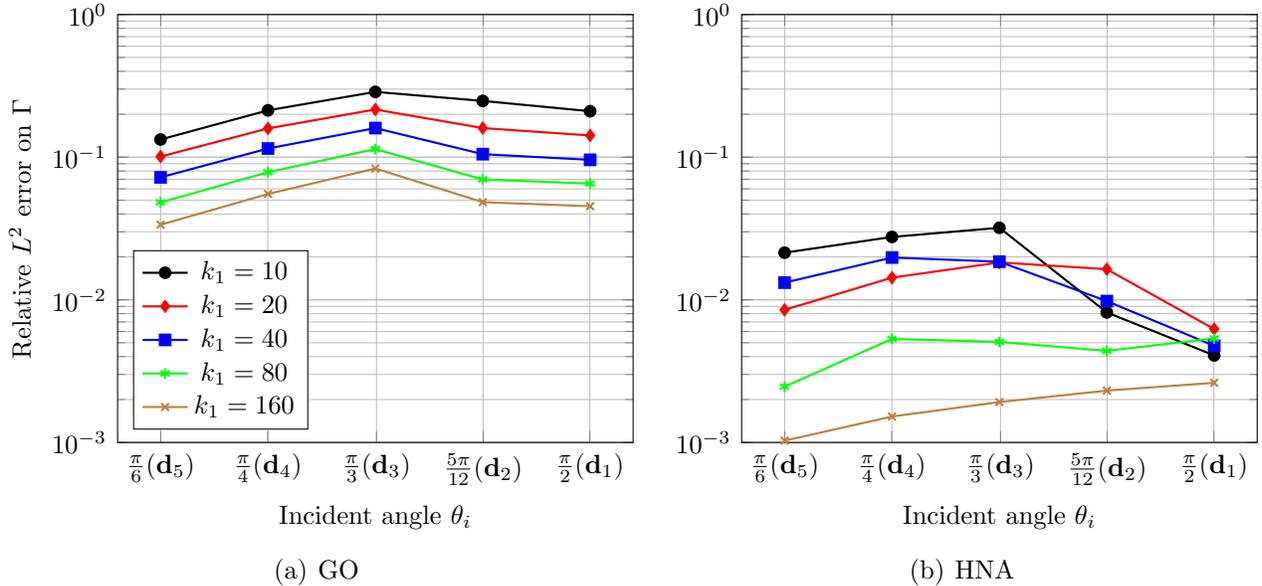
\begin{figure}[t!]
\centering
\subfigure[GO]{
\begin{tikzpicture}[scale=1]
\begin{semilogyaxis}[
	xlabel={Incident angle $\theta_i$},
	ylabel={Relative $L^2$ error on $\Gamma$},
	ymax={1*10^0},
	ymin={1*10^(-3)},
	grid=both,
	xtick={0.524,0.785,1.05,1.31,1.57},
	xticklabels={$\frac{\pi}{6}(\bd_5)$, $\frac{\pi}{4}(\bd_4)$, $\frac{\pi}{3}(\bd_3)$, $\frac{5\pi}{12}(\bd_2)$, $\frac{\pi}{2}(\bd_1)$},
	legend pos={south west},
	font=\footnotesize,
]

\addplot[color=black,mark=*] coordinates {
	(pi/6 , 1.33e-1)  (pi/4 , 2.13e-1)  (pi/3 , 2.87e-1)  (5*pi/12 , 2.48e-1)  (pi/2 , 2.10e-1)
};      
\addlegendentry{$k_1=10$}
	
\addplot[color=red,mark=diamond*] coordinates {
	(pi/6 , 1.01e-1)  (pi/4 , 1.59e-1)  (pi/3 , 2.16e-1)  (5*pi/12 , 1.60e-1)  (pi/2 , 1.42e-1)
};      
\addlegendentry{$k_1=20$}
	
\addplot[color=blue,mark=square*] coordinates {
	(pi/6 , 7.22e-2)  (pi/4 , 1.15e-1)  (pi/3 , 1.60e-1)  (5*pi/12 , 1.05e-1)  (pi/2 , 9.58e-2)
};      
\addlegendentry{$k_1=40$}
	
\addplot[color=green,mark=asterisk] coordinates {
	(pi/6 , 4.82e-2)  (pi/4 , 7.84e-2)  (pi/3 , 1.14e-1)  (5*pi/12 , 6.98e-2)  (pi/2 , 6.53e-2)
};   
\addlegendentry{$k_1=80$}
	
\addplot[color=brown,mark=x] coordinates {
	(pi/6 , 3.36e-2)  (pi/4 , 5.52e-2)  (pi/3 , 8.33e-2)  (5*pi/12 , 4.84e-2)  (pi/2 , 4.53e-2)
}; 
\addlegendentry{$k_1=160$}
	
\end{semilogyaxis}
\end{tikzpicture}
}
\subfigure[HNA]{
\begin{tikzpicture}[scale=1]
\begin{semilogyaxis}[
	xlabel={Incident angle $\theta_i$},
	ymax={1*10^0},
	ymin={1e-3},
	grid=both,
	legend columns={2},
	font=\footnotesize,
	xtick={0.524,0.785,1.05,1.31,1.57},
	xticklabels={$\frac{\pi}{6}(\bd_5)$, $\frac{\pi}{4}(\bd_4)$, $\frac{\pi}{3}(\bd_3)$, $\frac{5\pi}{12}(\bd_2)$, $\frac{\pi}{2}(\bd_1)$},
]

\addplot[color=black,mark=*] coordinates {
	(pi/6 , 2.14e-2)  (pi/4 , 2.76e-2)  (pi/3 , 3.20e-2)  (5*pi/12 , 8.16e-3)  (pi/2 , 4.08e-3)
};  

\addplot[color=red,mark=diamond*] coordinates {
	(pi/6 , 8.51e-3)  (pi/4 , 1.43e-2)  (pi/3 , 1.83e-2)  (5*pi/12 , 1.64e-2)  (pi/2 , 6.24e-3)
};   	
	
\addplot[color=blue,mark=square*] coordinates {
	(pi/6 , 1.32e-2)  (pi/4 , 1.98e-2)  (pi/3 , 1.85e-2)  (5*pi/12 , 9.80e-3)  (pi/2 , 4.76e-3)
};      
	
\addplot[color=green,mark=asterisk] coordinates {
	(pi/6 , 2.46e-3)  (pi/4 , 5.32e-3)  (pi/3 , 5.07e-3)  (5*pi/12 , 4.39e-3)  (pi/2 , 5.36e-3)
};   
	
\addplot[color=brown,mark=x] coordinates {
	(pi/6 , 1.03e-3)  (pi/4 , 1.52e-3)  (pi/3 , 1.92e-3)  (5*pi/12 , 2.31e-3)  (pi/2 , 2.62e-3)
};

\end{semilogyaxis}
\end{tikzpicture}
}

\caption{Accuracy of GO and HNA approximations to $u_1$ on $\Gamma$ plotted against incident angle $\theta_i$ for the incident directions $\bd_1,\ldots,\bd_5$ of Figure \ref{fig:triangle_setup}(a), for $\alpha=1$ and $\mu = 0.66+0.003125\ri$.
}
\label{fig:vary_angle_066}
\end{figure}

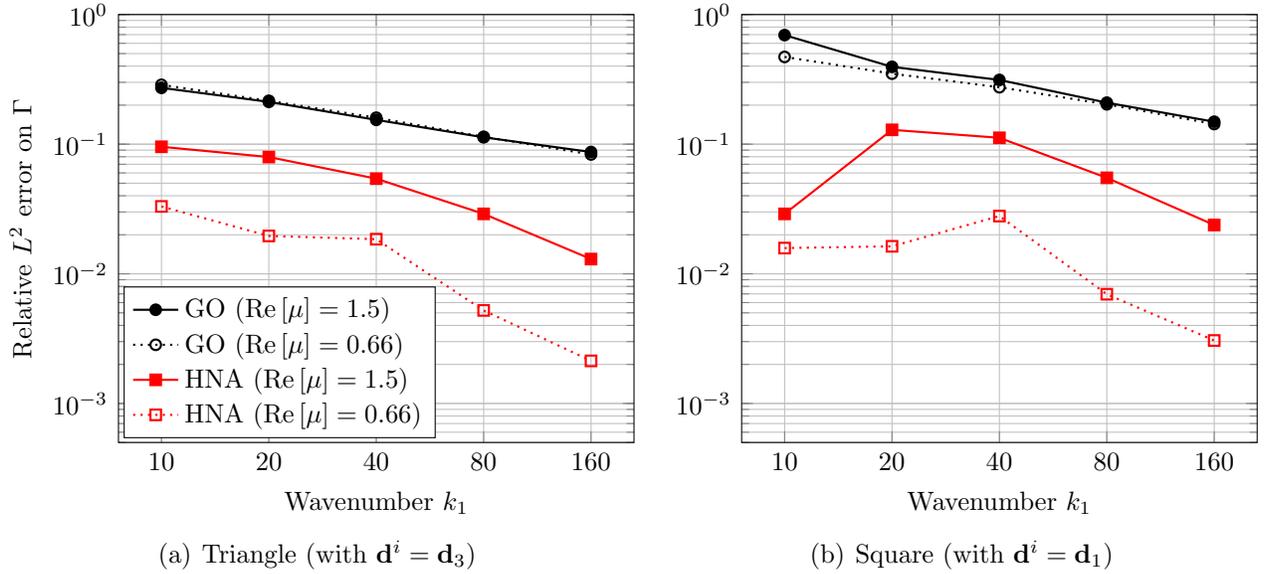
\begin{figure}[t!]
\centering
\subfigure[Triangle (with $\bd^i=\bd_3$)]{
\begin{tikzpicture}[scale=1]
\begin{loglogaxis}[
	xlabel={Wavenumber $k_1$},
	ylabel={Relative $L^2$ error on $\Gamma$},
	ymax={1e0},
	ymin={5e-4},
	grid=both,
	xtick={10,20,40,80,160},
	xticklabels={10,20,40,80,160},
	font=\footnotesize,
	legend cell align={left},
	legend style={at={(0.01,0.01)},anchor=south west}
]
	
\addplot[color=black,mark=*] coordinates {
	(10 , 2.72e-1)  (20 , 2.12e-1)  (40 , 1.54e-1)  (80 , 1.13e-1)  (160 , 8.70e-2) 
};
\addlegendentry{GO ($\real{\mu}=1.5$)}	

\addplot[color=black,dotted,mark=o,mark options={solid}] coordinates {
	(10 , 2.87e-1)  (20 , 2.16e-1)  (40 , 1.60e-1)  (80 , 1.14e-1)  (160 , 8.33e-2) 
};   
\addlegendentry{GO ($\real{\mu}=0.66$)}	

\addplot[color=red,mark=square*] coordinates {
	(10 , 9.55e-2)  (20 , 7.96e-2)  (40 , 5.42e-2)  (80 , 2.90e-2)  (160 , 1.30e-2) 
};
\addlegendentry{HNA ($\real{\mu}=1.5$)}
	
\addplot[color=red,dotted,mark=square,mark options={solid}] coordinates {
	(10 , 3.31e-2)  (20 , 1.96e-2)  (40 , 1.85e-2)  (80 , 5.22e-3)  (160 , 2.13e-3)  
};   
\addlegendentry{HNA ($\real{\mu}=0.66$)}	

\end{loglogaxis}
\end{tikzpicture}
}
\subfigure[Square (with $\bd^i=\bd_1$)]{
\begin{tikzpicture}[scale=1]
\begin{loglogaxis}[
	xlabel={Wavenumber $k_1$},
	ymax={1e0},
	ymin={5e-4},
	grid=both,
	font=\footnotesize,
	xtick={10,20,40,80,160},
	xticklabels={10,20,40,80,160},
	legend style={at={(0.01,0.01)},anchor=south west}
]
	
	\addplot[color=black,mark=*] coordinates {
	(10 , 6.94e-1)  (20 , 3.95e-1)  (40 , 3.13e-1)  (80 , 2.09e-1)  (160 , 1.49e-1)
	};
	
	\addplot[color=black,dotted,mark=o,mark options={solid}] coordinates {
	(10 , 4.71e-1)  (20 , 3.50e-1)  (40 , 2.75e-1)  (80 , 2.03e-1)  (160 , 1.43e-1)
	};   

	\addplot[color=red,mark=square*] coordinates {
	(10 , 2.90e-2)  (20 , 1.29e-1)  (40 , 1.12e-1)  (80 , 5.5e-2)  (160 , 2.38e-2)
	};
	
	\addplot[color=red,dotted,mark=square,mark options={solid}] coordinates {
	(10 , 1.58e-2)  (20 , 1.63e-2)  (40 , 2.79e-2)  (80 , 6.96e-3)  (160 , 3.06e-3)
	};   
\end{loglogaxis}
\end{tikzpicture}
}

\caption{Accuracy of GO and HNA approximations to $u_1$ on $\Gamma$ at grazing incidence, for (a) the triangle and (b) the square. We fix $\alpha=1$ and consider both $\mu = 1.5 + 0.003125\ri$ and $\mu = 0.66 + 0.003125\ri$. 
}
\label{fig:vary_contrast}
\end{figure}

In this section we study the performance of the HNA method in the case $0<\real{\mu}<1$. Specifically, fixing $\alpha=1$, we compare the results presented in Figure~\ref{fig:vary_angle} for the case $\mu = 1.5 + 0.003125\ri$, with those presented in Figure~\ref{fig:vary_angle_066} for $\mu = 0.66 + 0.003125\ri$ (the real part of the latter being approximately the reciprocal of the real part of the former). Based on these plots we make a number of observations. First, the accuracy of the GO approximation is roughly the same for the two refractive indices. Second, as was the case for $\real\mu=1.5$, the GO approximation for $\real\mu=0.66$ exhibits its largest errors at grazing incidence - again we ascribe this to a particularly prominent diffracted field along the grazing side in this case. Third, the HNA approximation consistently produces smaller errors for $\real\mu=0.66$ than for $\real\mu=1.5$. Fourth, the HNA error for $\real\mu=0.66$ is not obviously worse at grazing incidence than for other incident directions, in contrast to the $\real\mu=1.5$ case. To make a direct comparison easier, in Figure~\ref{fig:vary_contrast}(a) we plot the errors for the two cases ($\real\mu=1.5$ and $\real\mu=0.66$) at the grazing incidence angle $\theta_i=\pi/3$ ($\bd^i=\bd_3$), for a range of values of $k_1$. The GO errors are virtually identical for the two cases, but the HNA errors are significantly smaller for $\real\mu=0.66$ than for $\real\mu=1.5$.
These observations all support our hypothesis that the deterioration in accuracy of our HNA method near grazing incidence for $\real\mu>1$ is due to the presence of prominent lateral wave contributions, which are not present for $\real\mu<1$.

\begin{figure}[t!]
\centering
\subfigure[$\mu=1.5+0.003125\ri$]{%
       \includegraphics[trim=13cm 0 10cm 0,clip,width=0.45\textwidth]{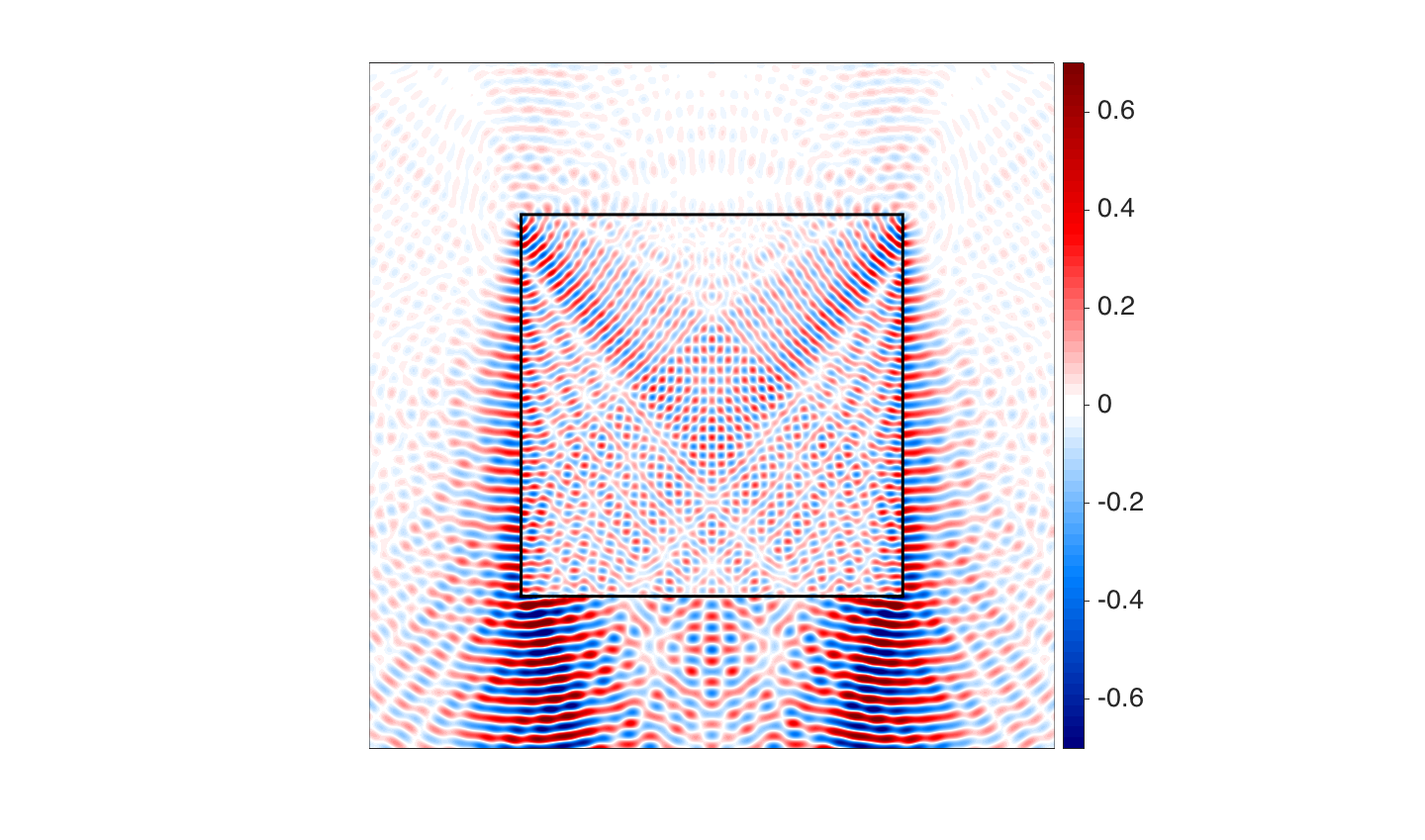}
    }
    \hspace{0.5cm}
    \subfigure[$\mu=0.66+0.003125\ri$]{%
       \includegraphics[trim=13cm 0 10cm 0,clip,width=0.45\textwidth]{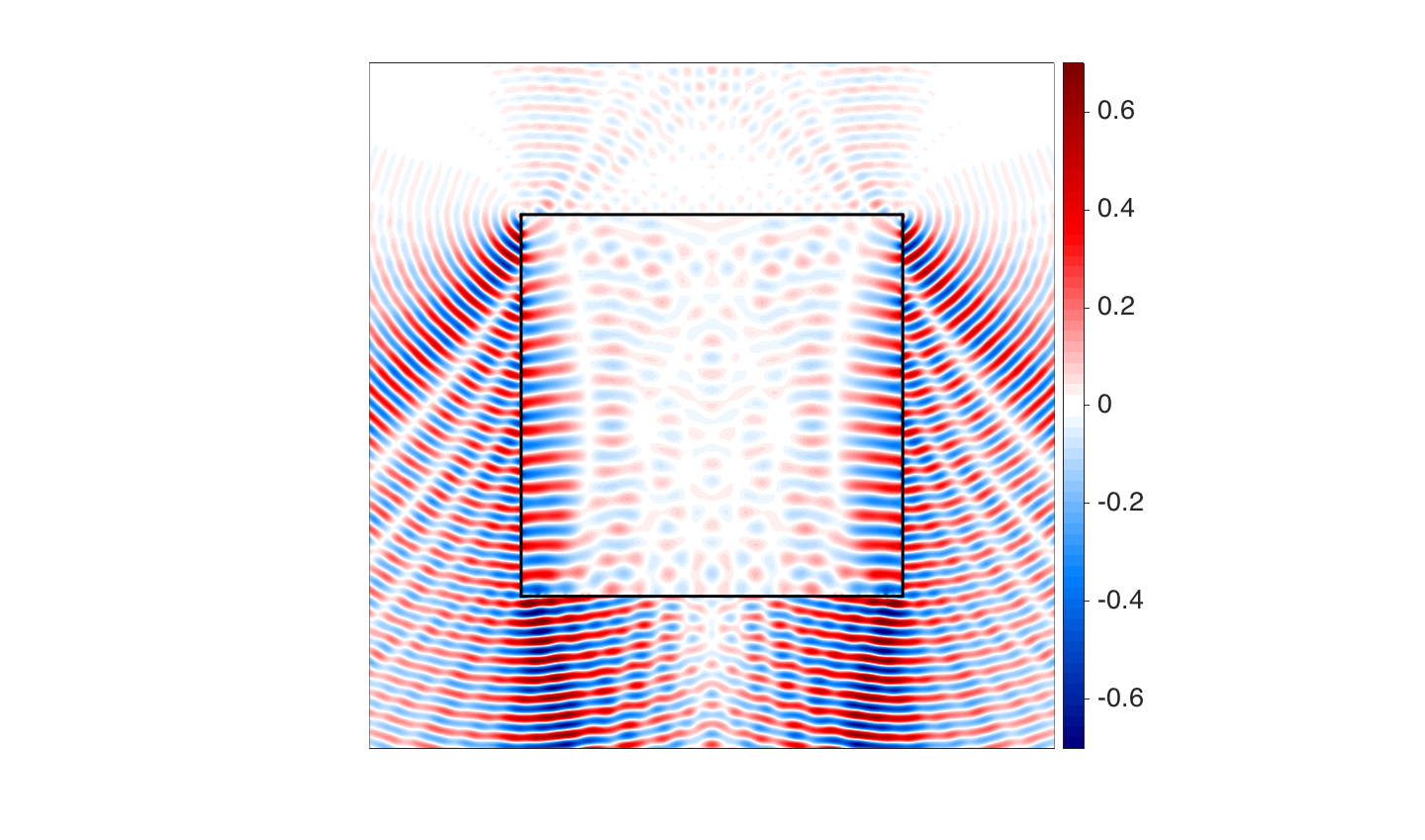}
    }
   \caption{Real part of the diffracted field for a square with different real part of the refractive index. The incident wave is directed from above.}
   \label{fig:sq_domain}
\end{figure}

To test this hypothesis further, we consider the analogous experiment for a square scatterer. Aligning the square with the Cartesian axes, the grazing incidence case is now $\bd^i=\bd_1$ ($\theta_i=\pi/2$). 
Field plots for $u_d$ in this configuration for the two cases ($\real\mu=1.5$ and $\real\mu=0.66$) can be found in Figure~\ref{fig:sq_domain}. 
The corresponding boundary errors are presented in Figure~\ref{fig:vary_contrast}(b). 
As for the triangle, the HNA error is significantly smaller for $\real\mu=0.66$ than for $\real\mu=1.5$, and again we hypothesize that this is due to the fact that for $\real\mu=1.5$ (and not for $\real\mu=0.66$) prominent lateral waves are generated by the two vertical sides, which impinge on the bottom side and limit the accuracy of our HNA approximation. 
We note also that this example allows us to rule out the possibility that increased error at grazing incidence might be due to the GO1/GO2 choice discussed in \S\ref{sec:comp_v_go}, because for this square configuration the propagation and decay directions ($\mathbf{d}$ and $\mathbf{e}$ in the notation of \S\ref{sec:comp_v_go}) always coincide and point vertically up or down, so that GO1 and GO2 coincide for all GO beams. 

\clearpage
\subsection{Varying $\alpha$}
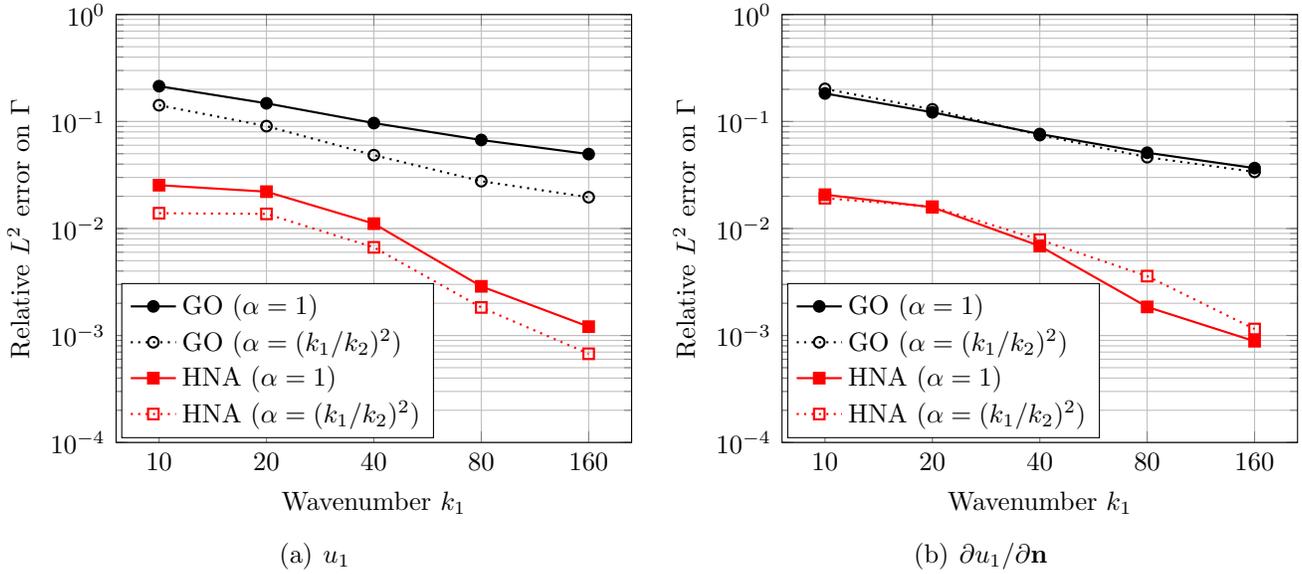
\begin{figure}
\centering
\subfigure[$u_1$]{
\begin{tikzpicture}[scale=1]
\begin{loglogaxis}[
	xlabel={Wavenumber $k_1$},
	ylabel={Relative $L^2$ error on $\Gamma$},
	ymax={1e0},
	ymin={1e-4},
	grid=both,
	font=\footnotesize,
	xtick={10,20,40,80,160},
	xticklabels={10,20,40,80,160},
	legend cell align={left},
	legend style={at={(0.01,0.01)},anchor=south west}
]
	
\addplot[color=black,mark=*] coordinates {
	(10 , 2.14e-1)  (20 , 1.48e-1)  (40 , 9.67e-2)  (80 , 6.72e-2)  (160 , 4.96e-2)  
};
\addlegendentry{GO ($\alpha=1$)}	

\addplot[color=black,dotted,mark=o,mark options={solid}] coordinates {
	(10 , 1.42e-1)  (20 , 9.06e-2)  (40 , 4.85e-2)  (80 , 2.77e-2)  (160 , 1.96e-2) 
};  
\addlegendentry{GO ($\alpha=(k_1/k_2)^2$)}	

\addplot[color=red,mark=square*] coordinates {
	(10 , 2.54e-2)  (20 , 2.21e-2)  (40 , 1.11e-2)  (80 , 2.88e-3)  (160 , 1.21e-3) 
};
\addlegendentry{HNA ($\alpha=1$)}
	
\addplot[color=red,dotted,mark=square,mark options={solid}] coordinates {
	(10 , 1.39e-2)  (20 , 1.37e-2)  (40 , 6.67e-3)  (80 , 1.83e-3)  (160 , 6.75e-4) 
};  
  
\addlegendentry{HNA ($\alpha=(k_1/k_2)^2$)}	
\end{loglogaxis}
\end{tikzpicture}
}
\subfigure[$\partial u_1/\partial\bn$]{
\begin{tikzpicture}[scale=1]
\begin{loglogaxis}[
	xlabel={Wavenumber $k_1$},
	ylabel={Relative $L^2$ error on $\Gamma$},
	ymax={1e0},
	ymin={1e-4},
	grid=both,
	font=\footnotesize,
	xtick={10,20,40,80,160},
	xticklabels={10,20,40,80,160},
	legend cell align={left},
	legend style={at={(0.01,0.01)},anchor=south west}
]
	
\addplot[color=black,mark=*] coordinates {
	(10 , 1.83e-1)  (20 , 1.22e-1)  (40 , 7.64e-2)  (80 , 5.11e-2)  (160 , 3.67e-2)  
};
\addlegendentry{GO ($\alpha=1$)}	

\addplot[color=black,dotted,mark=o,mark options={solid}] coordinates {
	(10 , 2.02e-1)  (20 , 1.30e-1)  (40 , 7.47e-2)  (80 , 4.64e-2)  (160 , 3.38e-2) 
};      
\addlegendentry{GO ($\alpha=(k_1/k_2)^2$)}	

\addplot[color=red,mark=square*] coordinates {
	(10 , 2.07e-2)  (20 , 1.58e-2)  (40 , 6.86e-3)  (80 , 1.85e-3)  (160 , 8.84e-4) 
};
\addlegendentry{HNA ($\alpha=1$)}
	
\addplot[color=red,dotted,mark=square,mark options={solid}] coordinates {
	(10 , 1.92e-2)  (20 , 1.60e-2)  (40 , 7.86e-3)  (80 , 3.59e-3)  (160 , 1.15e-3) 
};     
\addlegendentry{HNA ($\alpha=(k_1/k_2)^2$)}	
\end{loglogaxis}
\end{tikzpicture}
}
\caption{Accuracy of GO and HNA approximations to $u_1$ and $\pdonetext{u_1}{\bn}$ on $\Gamma$ for the triangle with $\alpha=1$ and $\alpha=(k_1/k_2)^2$, with $\mu=1.5+0.00625\ri$ and $\bd^i=\bd_1$ ($\theta_i=\pi/2$).}
\label{fig:vary_alpha}
\end{figure}

In principle, our HNA method can be evaluated for any value of the boundary condition jump parameter $\alpha\in\C\setminus \{0\}$. In practice, as was discussed after \rf{SRC}, the choices of $\alpha$ relevant to the dielectric electromagnetic scattering problem are 
$\alpha = 1$ and $\alpha = (k_1/k_2)^2$. We have already presented results for the former case; in this section we present results for the latter case. 
Specifically, we consider scattering by the triangle with $\mu=1.5+0.00625\ri$ and $\bd^i=\bd_1$ ($\theta_i=\pi/2$). 
Figure~\ref{fig:vary_alpha} shows that for this configuration the performance of the method in the approximation of both $u_1$ and $\partial u_1/\partial\bn$ is very similar for the two cases $\alpha = 1$ and $\alpha = (k_1/k_2)^2$. 

\subsection{Varying the absorption ($\im{\mu}$)}
\label{sec:vary_abs}
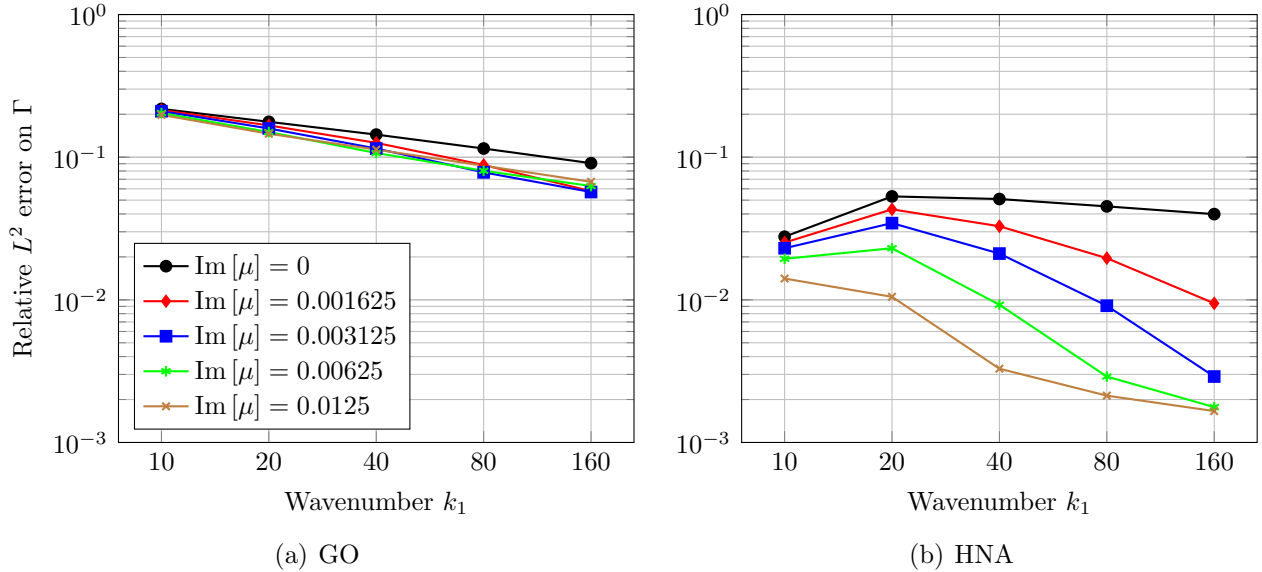
\begin{figure}[t!]
\centering
\subfigure[GO]{
\begin{tikzpicture}[scale=1]
\begin{loglogaxis}[
	xlabel={Wavenumber $k_1$},
	ylabel={Relative $L^2$ error on $\Gamma$},
	ymax={1*10^0},
	ymin={1e-3},
	grid=both,
	font=\footnotesize,
	xtick={10,20,40,80,160,320},
	xticklabels={10,20,40,80,160,320},
	legend pos={south west},
	legend cell align={left},
]
	
	\addplot[color=black,mark=*] coordinates {
	(10 , 2.18e-1)  (20 , 1.77e-1)  (40 , 1.44e-1)  (80 , 1.15e-1)  (160 , 9.08e-2) 
	}; 
	
	\addplot[color=red,mark=diamond*] coordinates {
	(10 , 2.14e-1)  (20 , 1.67e-1)  (40 , 1.26e-1)  (80 , 8.82e-2)  (160 , 5.78e-2) 
	};    
	
	\addplot[color=blue,mark=square*] coordinates {
	(10 , 2.10e-1)  (20 , 1.59e-1)  (40 , 1.15e-1)  (80 , 7.83e-2)  (160 , 5.70e-2) 
	}; 
	
	\addplot[color=green,mark=asterisk] coordinates {
	(10 , 2.04e-1)  (20 , 1.5e-1)  (40 , 1.07e-1)  (80 , 8.03e-2)  (160 , 6.28e-2) 
	};
	
	\addplot[color=brown,mark=x] coordinates {
	(10 , 1.98e-1)  (20 , 1.46e-1) (40 , 1.13e-1) (160 , 6.72e-2) 
	};
 
\legend{$\im\mu=0$,$\im\mu=0.001625$,$\im\mu=0.003125$,$\im\mu=0.00625$,$\im\mu=0.0125$}
\end{loglogaxis}
\end{tikzpicture}
}
\subfigure[HNA]{
\begin{tikzpicture}[scale=1]
\begin{loglogaxis}[
	xlabel={Wavenumber $k_1$},
	ymax={1e0},
	ymin={1e-3},
	grid=both,
	font=\footnotesize,
	xtick={10,20,40,80,160},
	xticklabels={10,20,40,80,160},
]
	
	\addplot[color=black,mark=*] coordinates {
	(10 , 2.77e-2)  (20 , 5.31e-2)  (40 , 5.09e-2)  (80 , 4.52e-2)  (160 , 3.99e-2)
	}; 
	
	\addplot[color=red,mark=diamond*] coordinates {
	(10 , 2.52e-2)  (20 , 4.31e-2)  (40 , 3.28e-2)  (80 , 1.96e-2)  (160 , 9.45e-3) 
	};    
	
	\addplot[color=blue,mark=square*] coordinates {
	(10 , 2.30e-2)  (20 , 3.45e-2)  (40 , 2.11e-2)  (80 , 9.09e-3)  (160 , 2.9e-3) 
	}; 
	
	\addplot[color=green,mark=asterisk] coordinates {
	(10 , 1.94e-2)  (20 , 2.30e-2)  (40 , 9.24e-3)  (80 , 2.90e-3)  (160 , 1.77e-3) 
	};
	
	\addplot[color=brown,mark=x] coordinates {
	(10 , 1.41e-2)  (20 , 1.05e-2)  (40 , 3.29e-3)  (80 , 2.13e-3)  (160 , 1.66e-3)
	};
 
\end{loglogaxis}
\end{tikzpicture}
}
\caption{Accuracy of GO and HNA approximations to $u_1$ on $\Gamma$ for the triangle with $\alpha=1$ and $\bd^i=\bd_5$ ($\theta_i=\pi/6$), with $\real\mu=1.5$ and $\im\mu$ varying between $0$ and $0.0125$.}
\label{fig:vary_abs}
\end{figure}

In this section we study the effect of varying the absorption $\im\mu$. 
Specifically, we consider scattering by the triangle with $\bd^i=\bd_5$ ($\theta_i=\pi/6$), $\alpha=1$, $\real{\mu}=1.5$ and 
\[
	\im{\mu}=0.0125,\,0.00625,\,0.003125,\,0.0015625,\, 0.
\]
In Figure~\ref{fig:vary_abs} 
we present the resulting error plots for the GO and HNA approximations. For all values of $\im\mu$ considered (even $\im\mu=0$) the HNA method provides an improvement in accuracy over the GO method and (once $k_1$ is sufficiently large) the HNA error decreases as $k_1$ increases. For fixed $k_1$ the HNA error decreases as the absorption $\im\mu$ increases; for the worst case ($\im\mu=0$) using the HNA approximation instead of GO alone reduces the error at $k_1=10$ from $21.2\%$ to $2.8\%$ and at $k_1=160$ from $9.1\%$ to $4.0\%$, whereas for the best case ($\im\mu=0.0125$) the error at $k_1=10$  is reduced from $19.8\%$ to $2.0\%$ and at $k_1=160$ from $6.7\%$ to $0.17\%$. 
This dependence of accuracy on absorption can be explained in terms of the asymptotic theory. Our HNA approximation space neglects the contribution of the higher-order asymptotic terms such as lateral waves and reflected-diffrated waves (see Figure \ref{fig:LateralAndDiffRef}), and for $\im\mu>0$ these higher-order terms are attenuated as they propagate across $\Omega_2$ at a rate proportional to $\im\mu$. Hence the larger $\im\mu$ is, the smaller the contribution the neglected terms make to the boundary solution $\bv$, and, as a result, the more accurate our HNA BEM is. 

\subsection{Scattering by a hexagon}
\label{subsec:geom}
\begin{figure}[t!]
\centering
\subfigure[Total field (real part)]{
      \includegraphics[trim=10cm 0 10cm 0,clip,width=0.5\textwidth]{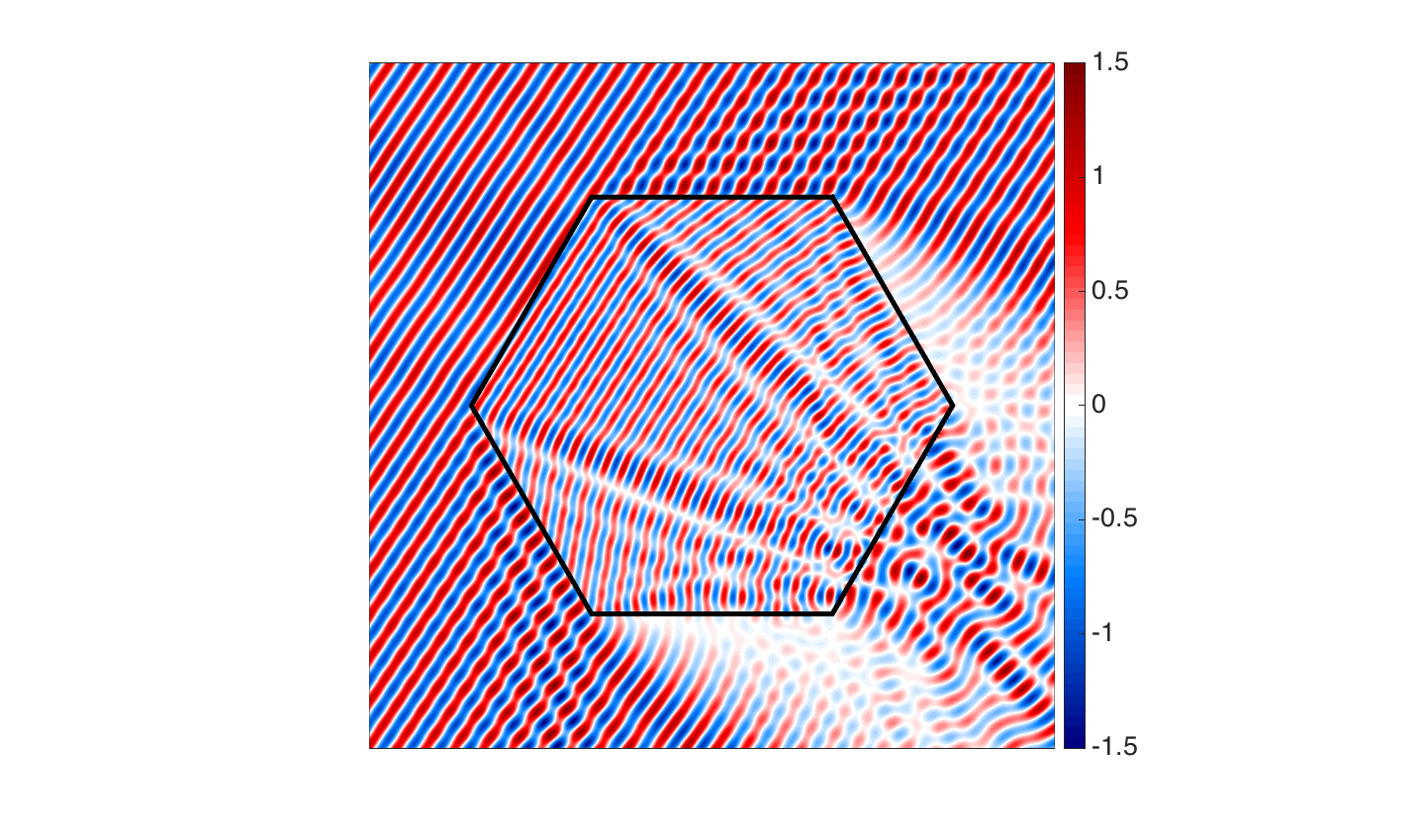}
   }
 \subfigure[Log of far-field pattern ($\log_{10}(|F|)$)]{
  \includegraphics[trim=27cm 0cm 25cm 1cm,clip,width=0.45\textwidth]{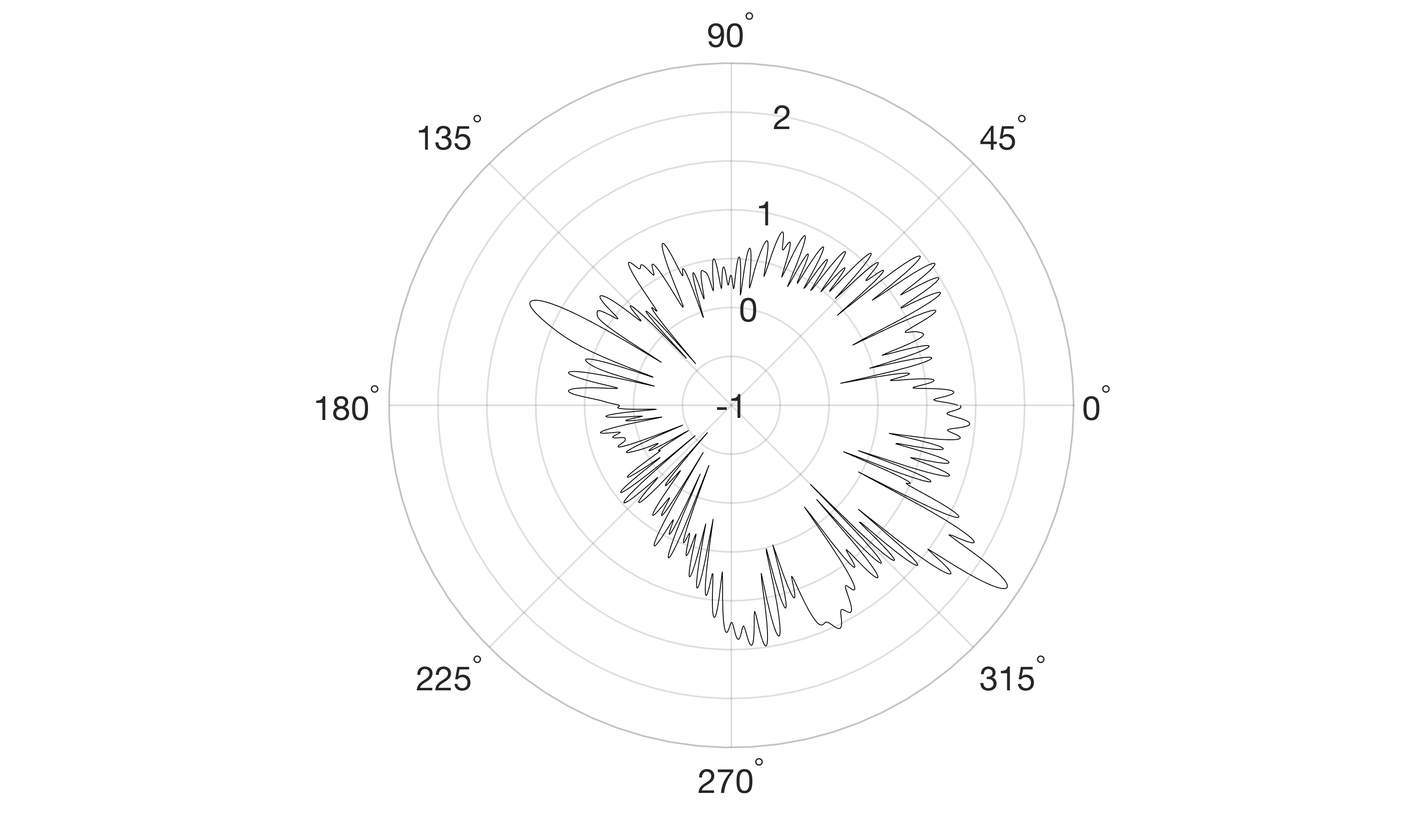}

 }
  \caption{Scattering by a regular hexagon with $k_1=10$, $\mu=1.39+0.00667\ri$, $\alpha=1$, and incident angle $\theta_i=\tan^{-1}(2/3)$.}
   \label{fig:hex_domain}
\end{figure}

We have already presented results for scattering by an equilateral triangle and a square. In this section we consider scattering by a regular hexagon. This can be viewed as a simple two-dimensional model of the motivating application we mentioned in \S\ref{sec:intro}, namely the scattering of light by atmospheric ice crystals, which typically take the form of hexagonal columns or their aggregates (see \cite{Baran:2009}). Here we choose a refractive index $\mu=1.39+0.00667\ri$ (which corresponds to the refractive index of ice for free-space wavelength 3.73$\mu m$ \cite{warren1984optical}), jump parameter $\alpha=1$ (corresponding to the out-of-plane electric field polarisation) and incident direction $\bd^i=(1/\sqrt{13})(3,-2)$ (incident angle $\theta_i = \tan^{-1}(2/3)$). The total field and far-field pattern for this configuration are plotted in Figure~\ref{fig:hex_domain}. 

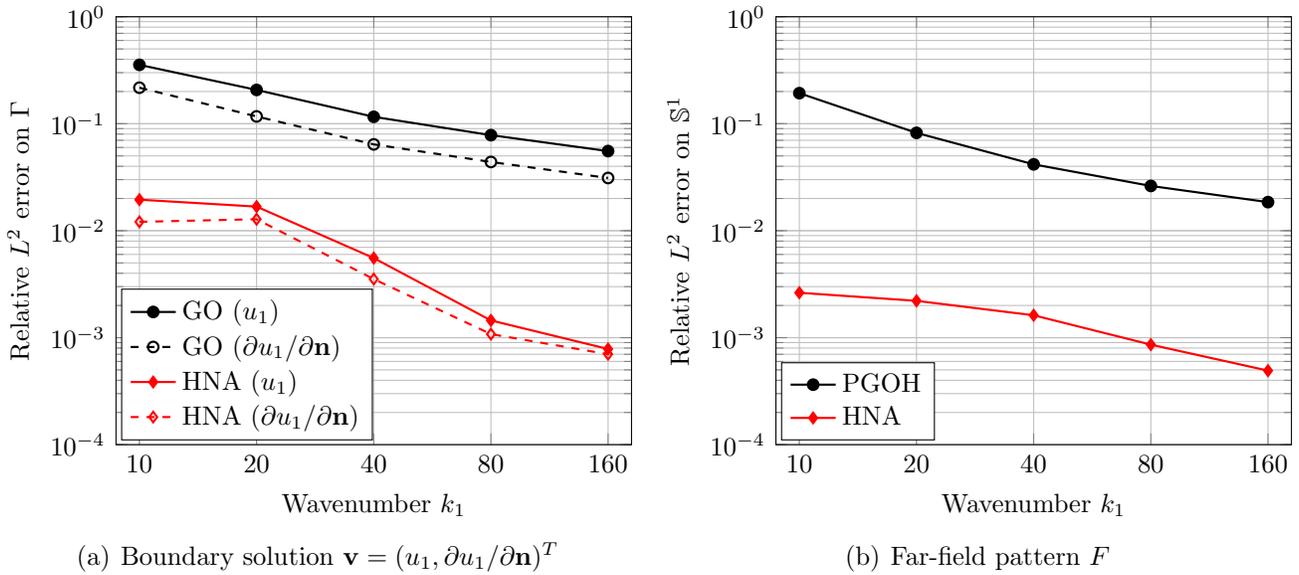
\begin{figure}[t!]
\centering
\subfigure[Boundary solution $\bv=(u_1,\pdonetext{u_1}{\bn})^T$]{
\begin{tikzpicture}[scale=1]
\begin{loglogaxis}[
	xlabel={Wavenumber $k_1$},
	ylabel={Relative $L^2$ error on $\Gamma$},
	ymax={1*10^0},
	ymin={1e-4},
	grid=both,
	font=\footnotesize,
	xtick={10,20,40,80,160},
	xticklabels={10,20,40,80,160},
	legend style={at={(0.01,0.01)},anchor=south west},
	legend cell align={left},
	enlarge x limits=0.05
]
\addplot[color=black,mark=*] coordinates {
	(10 , 3.55e-1)  (20 , 2.07e-1)  (40 , 1.16e-1)  (80 , 7.82e-2)  (160 , 5.55e-2)
};  
\addplot[color=black,dashed,mark=o,mark options={solid}]  coordinates {
	(10 , 2.17e-1)  (20 , 1.17e-1)  (40 , 6.41e-2)  (80 , 4.39e-2)  (160 , 3.11e-2)
}; 
	
\addplot[color=red,mark=diamond*] coordinates {
	(10 , 1.95e-2)  (20 , 1.68e-2)  (40 , 5.55e-3)  (80 , 1.45e-3)  (160 , 7.84e-4)
}; 
\addplot[color=red,dashed,mark=diamond,mark options={solid}] coordinates {
	(10 , 1.21e-2)  (20 , 1.28e-2)  (40 , 3.53e-3)  (80 , 1.08e-3)  (160 , 7.01e-4)
}; 
	 
\legend{GO ($u_1$),GO ($\partial u_1/\partial \bn$),HNA ($u_1$),HNA ($\partial u_1/\partial \bn$)}
\end{loglogaxis}
\end{tikzpicture}
}
\hs{-3}
\subfigure[Far-field pattern $F$]{
\begin{tikzpicture}[scale=1]
\begin{loglogaxis}[
	xlabel={Wavenumber $k_1$},
	ylabel={Relative $L^2$ error on $\mathbb{S}^1$},
	ymax={1*10^0},
	ymin={1e-4},
	grid=both,
	font=\footnotesize,
	xtick={10,20,40,80,160},
	xticklabels={10,20,40,80,160},
	legend style={at={(0.01,0.01)},anchor=south west},
	legend cell align={left},
	enlarge x limits=0.05
]
	
\addplot[color=black,mark=*] coordinates {
	(10 , 1.93e-1)  (20 , 8.21e-2)  (40 , 4.17e-2)  (80 , 2.62e-2)  (160 , 1.85e-2)
}; 
	
\addplot[color=red,mark=diamond*] coordinates {
	(10 , 2.63e-3)  (20 , 2.21e-3)  (40 , 1.62e-3)  (80 , 8.61e-4)  (160 , 4.92e-4)
};

\legend{PGOH,HNA}
\end{loglogaxis}
\end{tikzpicture}
}
\caption{Accuracy of GO and HNA for scattering by a regular hexagon with $\mu=1.39+0.00667\ri$, $\alpha=1$, and incident angle $\theta_i=\tan^{-1}(2/3)$.
}
\label{fig:hex_errors}
\end{figure}

The corresponding errors for $k_1=10,20,40,80,160$ are reported in Figure~\ref{fig:hex_errors}. The number of degrees of freedom $N$ in the HNA method decreases from 1008 at $k_1=10$ to 912 at $k_1=160$. (As was remarked in \S\ref{sec:vary_inc_angle}, our method typically uses \emph{fewer} degrees of freedom at larger wavenumbers because beam boundary discontinuities decrease in strength as the wavenumber increases.) 
For all wavenumbers considered, the HNA method provides a significant improvement in accuracy over GO (PGOH in the far field), and increases in accuracy as $k_1$ increases, despite the fact that the number of degrees of freedom $N$ remains fixed or decreases. Specifically, for the approximation of the far-field pattern at the highest wavenumber $k_1=160$, which corresponds to $960$ exterior wavelengths (approximately $1334$ interior wavelengths) around the scatterer boundary, the PGOH approximation achieves $1.9\%$ error, while the HNA BEM achieves $0.049\%$ error with just $N=912$, 
which corresponds to $\DOFtwo\approx 0.34$. 

\section{Discussion}
\label{sec:Discussion}
In this paper we have described an HNA BEM for high frequency scattering by 
penetrable (dielectric) convex polygons. We have demonstrated, by a range of numerical experiments, that our method can provide an order of magnitude improvement in accuracy over the purely asymptotic GO (PGOH in the far field) approach. Moreover, our method can achieve a fixed accuracy of approximation using a relatively small, frequency-independent number of degrees of freedom $N$. By contrast, conventional BEM approaches require $N$ to grow at least linearly with increasing frequency in order to maintain accuracy. The extremely efficient deployment of degrees of freedom in our HNA BEM is due to the fact that it uses an approximation space  built from oscillatory basis functions carefully chosen to capture certain components of the high frequency solution behaviour.

We conclude the paper with a discussion of some possible (but non-trivial) modifications that might enhance the performance of our method still further. First we recall that (as was explained in detail in \S\ref{sec:HNA}) our HNA method, while more accurate than GO alone, is not fully error-controllable (in the sense of the HNA methods analysed rigorously in \cite{CWL:07,hewett2013high,chandler2012high,hewett:shadow,hewett2014frequency}), because its accuracy is limited by our neglect of higher-order asymptotic effects such as lateral waves and diffracted-reflected waves (see Figure \ref{fig:LateralAndDiffRef}). Full HNA error-controllability seems infeasible for this problem because there are infinitely many different phases to consider, despite the fact that the scatterer is convex. But, in principle, the accuracy of our method could be improved, at the expense of additional degrees of freedom, by incorporating new oscillatory basis functions capable of capturing some of these higher-order effects. Concretely, for the lateral wave contributions one could include plane-wave basis functions of the form $f(\bx)\re^{\ri k_2 \bd^l \cdot \bx}$, where $f$ is a piecewise polynomial and $\bd^l$ is the propagation direction of the lateral wave in question (these directions can be determined solely from the refractive index and the corner angles of the polygon - for further discussion of these issues see \cite[\S3.2.1]{groth2015hybrid}). For the diffracted-reflected waves one could include basis functions of the form $f(\bx)\re^{\ri k_2 r_j'(\bx)}$, where $f$ is a piecewise polynomial and $r_j'(\bx)$ is the distance between $\bx$ and an appropriate reflected image $\bP_j'$ of the corner $\bP_j$ from which the original diffracted ray field emanated. (Such an image point $\bP_j'$ is shown in Figure \ref{fig:LateralAndDiffRef}(b).) We note that similar basis functions are required to capture multiple scattering effects in HNA BEMs for nonconvex impenetrable polygons - see, e.g., the discussion in \cite[\S8]{chandler2012high}.

Second, we suspect that our rather simplistic treatment of GO beam boundaries may not be optimal. It is well-known that, for scattering by impenetrable wedges, the smooth (but rapidly-varying) solution behaviour near the ``shadow boundaries'' across which GO components (incident or reflected) switch on/off is governed by the Fresnel integral and related functions (see, e.g., \cite{Ob:56,Kinber}). In the context of HNA methods for sound-soft nonconvex polygons this (frequency-dependent) shadow boundary behaviour can be captured either by including the appropriate special functions in the HNA approximation space (as in \cite{chandler2012high}) or by cutting off the GO component sharply across the shadow boundary and refining the mesh for the associated diffracted component towards it, to capture the rapid variation that compensates for the GO discontinuity (as in \cite{hewett:shadow}). 
Our approach in the current paper, in the context of penetrable convex polygons, is a rather crude implementation of the latter approach. We expect that certain minor (yet technical) modifications to our current practice might lead to slightly improved accuracy:
\begin{itemize}
\item Following \cite{groth2015hybrid}, we take the beam boundaries of a GO beam $a\re^{\ri (\NN\dr + \ri \KK\di)\cdot\bx}$ to be parallel to the associated propagation vector $\dr$. When $\KK>0$ (which will be the case generically if $\im\mu>0$) this choice may not be optimal - in simpler ``knife-edge'' diffraction problems involving inhomogeneous incident waves the ``correct'' location of shadow boundaries (defined in an appropriate sense) is known to be shifted in a rather non-intuitive way (see \cite{bertoni1978shadowing} and the discussion in \cite[Remark 3.1]{groth2015hybrid}). A possible modification to our current algorithm would be to allow an absorption-dependent shift in the beam boundary locations. 
\item For each beam boundary discontinuity classified as ``strong'' (in the sense of \rf{tolBBDef}), we insert a single new point in the mesh for the associated diffracted wave (if one exists). The rigorous analysis in \cite{hewett:shadow} for sound-soft nonconvex polygons suggests that the rapid solution variation near the beam boundary might be captured more accurately by introducing geometric mesh grading towards the beam boundary. 
\end{itemize}
To our knowledge, the canonical penetrable wedge scattering problem has not been analysed in sufficient detail to allow us to make theoretically-justifiable judgements about the optimal strategy in relation to the points raised above. Hence we leave such considerations
to future work.



\section*{Funding/acknowledgements}
This work was supported by EPSRC (EP/K000012/1 to SL, EP/P511262/1 to DPH, PhD studentship to SPG) and the Met Office (PhD CASE award to SPG). We would like to thank Anthony Baran for helpful comments and suggestions.

\bibliographystyle{siam}  
\bibliography{sams_bib}        
\end{document}

%% file: macros.tex
\newcommand{\done}[2]{\dfrac{d {#1}}{d {#2}}}
\newcommand{\donet}[2]{\frac{d {#1}}{d {#2}}}
\newcommand{\pdone}[2]{\dfrac{\partial {#1}}{\partial {#2}}}
\newcommand{\pdonet}[2]{\frac{\partial {#1}}{\partial {#2}}}
\newcommand{\pdonetext}[2]{\partial {#1}/\partial {#2}}
\newcommand{\pdtwo}[2]{\dfrac{\partial^2 {#1}}{\partial {#2}^2}}
\newcommand{\pdtwot}[2]{\frac{\partial^2 {#1}}{\partial {#2}^2}}
\newcommand{\pdtwomix}[3]{\dfrac{\partial^2 {#1}}{\partial {#2}\partial {#3}}}
\newcommand{\pdtwomixt}[3]{\frac{\partial^2 {#1}}{\partial {#2}\partial {#3}}}
\newcommand{\bs}[1]{\mathbf{#1}}
\newcommand{\bx}{\mathbf{x}}
\newcommand{\by}{\mathbf{y}}
\newcommand{\bd}{\mathbf{d}} 
\newcommand{\bn}{\mathbf{n}} 
\newcommand{\bP}{\mathbf{P}} 
\newcommand{\bp}{\mathbf{p}} 
\newcommand{\ol}[1]{\overline{#1}}
\newcommand{\rf}[1]{(\ref{#1})}
\newcommand{\xt}{\mathbf{x},t}
\newcommand{\hs}[1]{\hspace{#1mm}}
\newcommand{\vs}[1]{\vspace{#1mm}}
\newcommand{\eps}{\varepsilon}
\newcommand{\ord}[1]{\mathcal{O}\left(#1\right)} 
\newcommand{\oord}[1]{o\left(#1\right)}
\newcommand{\Ord}[1]{\Theta\left(#1\right)}
\newcommand{\PhiF}{\Phi_{\rm freq}}
\newcommand{\real}[1]{{\rm Re}\left[#1\right]} 
\newcommand{\im}[1]{{\rm Im}\left[#1\right]}
\newcommand{\hsnorm}[1]{||#1||_{H^{s}(\bs{R})}}
\newcommand{\hnorm}[1]{||#1||_{\tilde{H}^{-1/2}((0,1))}}
\newcommand{\norm}[2]{\left\|#1\right\|_{#2}}
\newcommand{\normt}[2]{\|#1\|_{#2}}
\newcommand{\on}[1]{\Vert{#1} \Vert_{1}}
\newcommand{\tn}[1]{\Vert{#1} \Vert_{2}}
\newcommand{\ts}{\tilde{s}}
\newcommand{\tGamma}{{\tilde{\Gamma}}}
\newcommand{\darg}[1]{\left|{\rm arg}\left[ #1 \right]\right|}
\newcommand{\bnabla}{\boldsymbol{\nabla}}
\newcommand{\dive}{\boldsymbol{\nabla}\cdot}
\newcommand{\curl}{\boldsymbol{\nabla}\times}
\newcommand{\Phixy}{\Phi(\bx,\by)}
\newcommand{\PhiOxy}{\Phi_0(\bx,\by)}
\newcommand{\dxPhixy}{\pdone{\Phi}{n(\bx)}(\bx,\by)}
\newcommand{\dyPhixy}{\pdone{\Phi}{n(\by)}(\bx,\by)}
\newcommand{\dxPhiOxy}{\pdone{\Phi_0}{n(\bx)}(\bx,\by)}
\newcommand{\dyPhiOxy}{\pdone{\Phi_0}{n(\by)}(\bx,\by)}

\newcommand{\rd}{\mathrm{d}}
\newcommand{\R}{\mathbb{R}}
\newcommand{\N}{\mathbb{N}}
\newcommand{\Z}{\mathbb{Z}}
\newcommand{\C}{\mathbb{C}}
\newcommand{\K}{{\mathbb{K}}}
\newcommand{\ri}{{\mathrm{i}}}
\newcommand{\re}{{\mathrm{e}}} 

\newcommand{\cA}{\mathcal{A}}
\newcommand{\cC}{\mathcal{C}}
\newcommand{\cS}{\mathcal{S}}
\newcommand{\cD}{\mathcal{D}}
\newcommand{\cone}{{c_{j}^\pm}}
\newcommand{\ctwo}{{c_{2,j}^\pm}}
\newcommand{\cthree}{{c_{3,j}^\pm}}

\newtheorem{thm}{Theorem}[section]
\newtheorem{lem}[thm]{Lemma}
\newtheorem{defn}[thm]{Definition}
\newtheorem{prop}[thm]{Proposition}
\newtheorem{cor}[thm]{Corollary}
\newtheorem{rem}[thm]{Remark}
\newtheorem{conj}[thm]{Conjecture}
\newtheorem{ass}[thm]{Assumption}
\newtheorem{example}[thm]{Example} 